\documentclass[12pt,a4paper]{article}

\usepackage{graphicx}
\usepackage{amsfonts}
\usepackage{amsmath}
\usepackage{bigstrut}

\makeatletter
\newdimen\mynormalparindent
\def\mymakefnmark{}
\def\mymakefntext{\indent\mymakefnmark}
\long\def\myfootnotetext#1{\insert\footins{%
  \normalfont\footnotesize
  \interlinepenalty\interfootnotelinepenalty
  \splittopskip\footnotesep \splitmaxdepth \dp\strutbox
  \floatingpenalty\@MM \hsize\columnwidth
  \@parboxrestore \parindent\mynormalparindent \sloppy
  \mymakefntext{\rule\z@\footnotesep\ignorespaces#1\unskip\strut\par}}}

\makeatother

\long\def\symbolfootnote[#1]#2{\begingroup%
\def\thefootnote{\fnsymbol{footnote}}\footnote[#1]{#2}\endgroup}

\newtheorem{basic}{Basic}[section]

\newtheorem{lem}[basic]{Lemma}
\newtheorem{propos}[basic]{Proposition}
\newtheorem{thm}[basic]{Theorem}

\newtheorem{cor}[basic]{Corollary}

\newcommand{\bdm}{\begin{displaymath}}
\newcommand{\edm}{\end{displaymath}}
\newcommand{\be}{\begin{equation}}
\newcommand{\ee}{\end{equation}}
\newcommand{\ep}{\vspace{-3mm}\hfill\mbox{$\Box$}\\}

\newcommand{\ol}{\overline}
\newcommand{\R}{\mathbb{R}}
\newcommand{\C}{\mathbb{C}}
\newcommand{\wt}{\widetilde}
\newcommand{\wh}{\widehat}

\def\theequation{\thesection.\arabic{equation}}

\begin{document}

\begin{center}
{\bf\large On best rank-2 and rank-(2,2,2) approximations of
order-3 tensors}\\

\vspace{2cm} Alwin Stegeman$^*$
\hspace{5mm} and\hspace{5mm} 
Shmuel Friedland$^{\dagger}$
\end{center}
\vspace{2cm}

$^*$ Corresponding author, Heijmans Institute for
Psychological Research, University of Groningen, Grote Kruisstraat 2/1,
9712 TS Groningen, The Netherlands, phone: ++31 597 551308, The Netherlands, {\tt stegeman.alwin@gmail.com}, {\tt http://www.alwinstegeman.nl}.\\~\vspace{5mm}

$^\dagger$ Department of Mathematics, Statistics and Computer Science, University of Illinois at Chicago, Chicago, Illinois 60607-7045, USA, \texttt{friedlan@uic.edu}. \\~\vspace{4cm}

\noindent {\bf Acknowledgement.} This work was supported in part by the National Science Foundation under Grant DMS-1216393.

\newpage
\begin{center}
{\bf\large On best rank-2 and rank-(2,2,2) approximations of
order-3 tensors}\\
\vspace{1cm}
\begin{abstract}
\noindent It is well known that a best rank-$R$ approximation of order-3 tensors may not exist for $R\ge 2$. A best rank-$(R,R,R)$ approximation always exists, however, and is also a best rank-$R$ approximation when it has rank (at most) $R$. For $R=2$ and real order-3 tensors it is shown that a best rank-2 approximation is also a local minimum of the best rank-(2,2,2) approximation problem. This implies that if all rank-(2,2,2) minima have rank larger than 2, then a best rank-2 approximation does not exist. This provides an easy-to-check criterion for existence of a best rank-2 approximation. The result is illustrated by means of simulations. \\~\\

\end{abstract}\end{center}

\section{Introduction}
\setcounter{equation}{0}
Order-$N$ tensors are defined on the Cartesian product of $N$ linear spaces. For fixed bases of these linear spaces the tensor is represented by an $N$-way array. Tensors are usually identified by their array representation. In this paper we consider real order-3 tensors. The rank of an order-3 tensor ${\cal Y}\in\R^{I\times J\times K}$ is defined as 
\be
{\rm rank}({\cal Y})=\min\{R\,:\,{\cal Y}=\sum_{r=1}^R ({\bf a}_r\circ{\bf b}_r\circ{\bf c}_r)\}\,,
\ee

\noindent where ${\bf a}_r\in\R^I$, ${\bf b}_r\in\R^J$, ${\bf c}_r\in\R^K$, $r=1,\ldots,R$, and $\circ$ denotes the outer vector product. The outer vector product ${\cal Y}={\bf a}\circ{\bf b}\circ{\bf c}$ has entries $y_{ijk}=a_i\,b_j\,c_k$ and constitutes a rank-1 tensor when ${\bf a}$, ${\bf b}$, and ${\bf c}$ are nonzero.  The set of tensors with rank at most $R$ is denoted by
\be
S_R(I,J,K)=\{{\cal Y}\in\R^{I\times J\times K}:\;{\rm rank}({\cal Y})\le R\}\,.
\ee

\noindent We consider the problem of finding a best rank-$R$ approximation of ${\cal Z}\in\R^{I\times J\times K}$:
\be
\label{prob-CPD}
\min_{{\cal Y}\in S_R(I,J,K)} \|{\cal Z}-{\cal Y}\|^2_F\,,
\ee

\noindent where $\|\cdot\|_F$ denotes the Frobenius norm (i.e., the square root of the sum-of-squares). The problem (\ref{prob-CPD}) is solved in the variables ${\bf a}_r\in\R^I$, ${\bf b}_r\in\R^J$, ${\bf c}_r\in\R^K$, $r=1,\ldots,R$, and its solution is known as a Canonical Polyadic Decomposition (CPD), or Candecomp \cite{CC}, or Parafac \cite{Har}. It was introduced in \cite{H1,H2}. The CPD and related decompositions have numerous applications \cite{SBG,Kro,CoLa,DL5,Bini-1,Bini-2,KB,Acar}, and various iterative CPD algorithms are available \cite{TB}. Unfortunately, for $R\ge 2$ the problem may not have an optimal solution because the set $S_R(I,J,K)$ is not closed \cite{DSL}. In such a case, trying to compute a best rank-$R$ approximation yields a rank-$R$ sequence converging to a boundary point ${\cal X}$ of $S_R(I,J,K)$ with rank$({\cal X})>R$. This is accompanied by rank-1 terms that become nearly linearly dependent, while their norms become arbitrarily large \cite{KHL,KDS,DSL}. This phenomenon is known as ``diverging components" or ``degenerate solutions" or ``diverging rank-1 terms" \cite{Paa2,Ste,Ste3,SDL-TR}. To guarantee existence of a best rank-$R$ decomposition, one may impose constraints on the rank-1 terms such as orthogonality or nonnegativity \cite{KDS,LiCo,LiCo2}. However, these constraints are not suitable for all applications.

Results on nonexistence of a best rank-$R$ approximation are the following. Any $2\times 2\times 2$ tensor of rank 3 does not have a best rank-2 approximation \cite{DSL}. Nonexistence of a best rank-2 approximation of ${\cal Z}\in\R^{I\times J\times K}$ holds on a set of positive volume \cite{DSL}. Nonexistence of a best rank-$R$ approximation holds on a set of positive volume or even almost everywhere for certain classes of ${\cal Z}\in\R^{I\times J\times 2}$ \cite{Ste2,Ste-arxiv-IJ2}. Note that a best rank-1 approximation always exists, since $S_1(I,J,K)$ is closed \cite{DSL}.

Instead of imposing constraints, one may consider the following problem instead:
\be
\label{prob-CPDc}
\min_{{\cal Y}\in \ol{S}_R(I,J,K)} \|{\cal Z}-{\cal Y}\|^2_F\,,
\ee

\noindent where $\ol{S}_R(I,J,K)$ denotes the closure of $S_R(I,J,K)$ in $\R^{I\times J\times K}$. To solve problem (\ref{prob-CPDc}), a characterization of the boundary points of $S_R(I,J,K)$ is needed, and an algorithm to find an optimal boundary point. For $S_R(I,J,2)$ and $R\le\min(I,J)$ this can be done via the Generalized Schur Decomposition (GSD) \cite{SDL,Ste-arxiv-GSD}. Results on the existence of best rank-$R$ approximations for generic ${\cal Z}\in\R^{I\times J\times 2}$ can be found in \cite{Ste2,Ste-arxiv-IJ2}. For $S_2(I,J,K)$ the boundary points are described in \cite{DSL} and an algorithm is developed in \cite{RoG} by means of finding a best rank-(2,2,2) approximation with several zero restrictions on the $2\times 2\times 2$ core tensor. If the free (2,2,2)-entry of the core tensor equals zero for all such best rank-(2,2,2) approximations, then the optimal boundary point of $\ol{S}_2(I,J,K)$ has rank 3 and no best rank-2 approximation exists. 

A different approach to find an optimal boundary point of $\ol{S}_R(I,J,K)$ is to determine its decomposition form from the pattern of groups of diverging rank-1 terms in the CPD sequence generated by the iterative CPD algorithm. The decomposition can then be fitted to the tensor using initial values obtained from the diverging CPD sequence. This method has been proposed and demonstrated in simulation studies in \cite{Ste-limit,Ste-Jordan}. For an application see \cite{Ste-TV}.

In this paper, we prove a new criterion for existence of a best rank-2 approximation for real order-3 tensors via (unconstrained) best rank-(2,2,2) approximations as recently suggested in \cite{FriTam15}. The $2\times 2\times 2$ core tensor of a best rank-(2,2,2) approximation has rank 2 or 3. When it has rank 2 the best rank-(2,2,2) approximation is also a best rank-2 approximation. We prove that a best rank-2 approximation is also a local minimum of the best rank-(2,2,2) approximation problem for generic tensors. Hence, if all best rank-(2,2,2) minima have rank 3, then no best rank-2 approximation exists. Verifying the rank of the $2\times 2\times 2$ core array is numerically more reliable than checking whether a core entry equals zero \cite{RoG}. This will be demonstrated in a simulation study. A data analytic perspective on the equivalence of a best rank-(2,2,2) approximation of rank 2 and a best rank-2 approximation can be found in \cite{KTB}.

The paper is organized as follows. In section 2 we consider the problem of finding a best rank-$(R_1,R_2,R_3)$ approximation by iterating over Givens rotations to find the orthonormal bases for the three subspaces. The first-order conditions of this problem are derived. In section 3 we parameterize the set $\ol{S}_2(I,J,K)$ by a variant of the GSD and derive first-order conditions using Givens rotations as in section 2. The first-order conditions are used to prove our main result in section 4. In section 5 we compare our criterion for existence of a best rank-2 approximation to that of \cite{RoG} in a simulation study. In section 6 we consider the case of complex tensors and provide links between our results and existing results in algebraic geometry. Finally, section 7 contains a discussion of our findings.

We use the following notation. The notation ${\cal Y}$, ${\bf Y}$, ${\bf
y}$, $y$ is used for a three-way array, a matrix, a column vector, and a
scalar, respectively. All arrays, matrices, vectors, and scalars are
real-valued unless indicated otherwise. Matrix transpose and inverse are denoted as ${\bf
Y}^T$ and ${\bf Y}^{-1}$, respectively. A zero matrix of size $p\times
q$ is denoted by ${\bf O}_{p,q}$. A zero column vector is
denoted by ${\bf 0}$. A $p\times p$ matrix ${\bf Y}$ is called orthogonal 
if ${\bf Y}^T{\bf Y}={\bf YY}^T={\bf I}_p$. A $p\times q$ matrix, $p>q$, is called columnwise orthogonal if ${\bf Y}^T{\bf Y}={\bf I}_q$.

\section{Finding a best rank-$(R_1,R_2,R_3)$ approximation}
\setcounter{equation}{0}
A mode-$i$ vector or fiber of an order-3 tensor ${\cal Y}\in\R^{I\times J\times K}$ is given by varying the $i$th index while keeping the other two indices fixed. The mode-$i$ rank, denoted by rank$_i({\cal Y})$, is defined as the rank of the collection of mode-$i$ vectors. Hence, for generic ${\cal Y}\in\R^{I\times J\times K}$ we have rank$_1({\cal Y})=\min(I,JK)$, rank$_2({\cal Y})=\min(J,IK)$, and rank$_3({\cal Y})=\min(K,IJ)$. The triplet $({\rm rank}_1({\cal Y}),{\rm rank}_2({\cal Y}),{\rm rank}_3({\cal Y}))$ is referred to as the multilinear rank, which we denote as mrank$({\cal Y})$. The set of tensors of multilinear rank at most $(R_1,R_2,R_3)$ is denoted by
\be
M_{(R_1,R_2,R_3)}(I,J,K)=\{{\cal Y}\in\R^{I\times J\times K}:\;{\rm mrank}({\cal Y})\le (R_1,R_2,R_3)\}\,,
\ee

\noindent where we assume $R_1\le I$, $R_2\le J$, and $R_3\le K$. We define the multilinear transformation ${\cal Y}=({\bf S},{\bf T},{\bf U})\cdot{\cal G}$ via $y_{ijk}=\sum_{pqr}s_{ip}t_{jq}u_{kr}g_{pqr}$. Any ${\cal Y}\in M_{(R_1,R_2,R_3)}(I,J,K)$ can be written as ${\cal Y}=({\bf S},{\bf T},{\bf U})\cdot{\cal G}$, with ${\bf S}\in\R^{I\times R_1}$, ${\bf T}\in\R^{J\times R_2}$, and ${\bf U}\in\R^{K\times R_3}$ being columnwise orthogonal, and ${\cal G}\in\R^{R_1\times R_2\times R_3}$. We consider the problem of finding a best rank-$(R_1,R_2,R_3)$ approximation to a given tensor ${\cal Z}\in\R^{I\times J\times K}$:
\be
\label{prob-Mult}
\min_{{\cal Y}\in M_{(R_1,R_2,R_3)}(I,J,K)} \|{\cal Z}-{\cal Y}\|^2_F\,,
\ee

\noindent which is solved in the variables ${\bf S}$, ${\bf T}$, ${\bf U}$, and ${\cal G}$. A solution to problem (\ref{prob-Mult}) is also known as a Tucker3 decomposition \cite{Tuc}, with the higher-order singular value decomposition (HOSVD) \cite{LMV-hosvd} being a Tucker3 solution in which the transformational ambiguities are (mostly) fixed.

Problem (\ref{prob-Mult}) is equivalent to maximizing $\|({\bf S}^T,{\bf T}^T,{\bf U}^T)\cdot{\cal Z}\|^2_F$ over $({\bf S},{\bf T},{\bf U})$ and setting ${\cal G}=({\bf S}^T,{\bf T}^T,{\bf U}^T)\cdot{\cal Z}$; see \cite{LMV00}. This in turn is equivalent to maximizing the Frobenius norm of any $R_1\times R_2\times R_3$ subtensor of $(\wt{\bf S}^T,\wt{\bf T}^T,\wt{\bf U}^T)\cdot{\cal Z}$ over orthogonal $\wt{\bf S}\in\R^{I\times I}$, $\wt{\bf T}\in\R^{J\times J}$, and $\wt{\bf U}\in\R^{K\times K}$, with ${\cal G}$ equal to the transformed subtensor, and $({\bf S},{\bf T},{\bf U})$ equal to the corresponding columns of $(\wt{\bf S},\wt{\bf T},\wt{\bf U})$. We work with this formulation of problem (\ref{prob-Mult}), where we take the subtensor that has indices $(i,j,k)$ with $1\le i\le R_1$, $1\le j\le R_2$, and $1\le k\le R_3$. Since the set of orthogonal matrices is compact, problem (\ref{prob-Mult}) is guaranteed to have an optimal solution.

Algorithms for solving problem (\ref{prob-Mult}) have been proposed in \cite{KrL,LMV00,Sav,Ish,FriTam15}. In the higher-order power method of \cite{LMV00} each iteration rotates mass to the target subtensor by using the singular value decomposition (SVD) of the columns of one of the three matrix unfoldings of the rotated tensor corresponding to the subtensor. First-order conditions are derived in \cite{ES} and correspond to orthogonal sets of vectors in the three matrix unfoldings. Let the slices of $\wt{\cal Z}=(\wt{\bf S}^T,\wt{\bf T}^T,\wt{\bf U}^T)\cdot{\cal Z}$ be partitioned as $\left[\begin{array}{cc}{\bf G}_k & {\bf L}_k\\ {\bf N}_k & {\bf M}_k\end{array}\right]$, where ${\bf G}_k\in\R^{R_1\times R_2}$, ${\bf L}_k\in\R^{R_1\times (J-R_2)}$, ${\bf N}_k\in\R^{(I-R_1)\times R_2}$, and ${\bf M}_k\in\R^{(I-R_1)\times (J-R_2)}$, $k=1,\ldots,K$. Hence, in an optimal solution the core tensor ${\cal G}$ has slices ${\bf G}_k$, $k=1,\ldots,R_3$. The first-order conditions can be written as
\be
\label{eq-statM-row}
[{\bf G}_1\;\cdots\;{\bf G}_{R_3}]\;[{\bf N}_1\;\cdots\;{\bf N}_{R_3}]^T={\bf O}\,,
\ee
\be
\label{eq-statM-col}
\left[\begin{array}{c}{\bf G}_1\\ \vdots \\ {\bf G}_{R_3}\end{array}\right]^T\;
\left[\begin{array}{c}{\bf L}_1\\ \vdots \\ {\bf L}_{R_3}\end{array}\right]={\bf O}\,,
\ee
\be
\label{eq-statM-slice}
[{\rm vec}({\bf G}_1)\;\ldots\;{\rm vec}({\bf G}_{R_3})]^T\;
[{\rm vec}({\bf G}_{R_3+1})\;\ldots\;{\rm vec}({\bf G}_K)]={\bf O}\,,
\ee

\noindent where vec$(\cdot)$ stacks the columns of a matrix below each other in a column vector. 

For later use, we present an alternative derivation of the first-order conditions. We consider updating each of $(\wt{\bf S},\wt{\bf T},\wt{\bf U})$ by means of Givens rotations, an approach that has been used for the Simultaneous Generalized Schur Decomposition (SGSD) in \cite{LMV}, and for the GSD in \cite{Ste-arxiv-IJ2}. The first-order conditions can be derived by requiring that rotating two rows, columns, or slices of $(\wt{\bf S}^T,\wt{\bf T}^T,\wt{\bf U}^T)\cdot{\cal Z}$ will not increase the Frobenius norm of its $R_1\times R_2\times R_3$ subtensor. We need the following result.

\begin{lem}
\label{lem-statorth}
For vectors ${\bf x},{\bf y}\in\R^p$ and $\alpha\in\R$, define the rotation
\bdm
[\tilde{\bf x}\;\;\tilde{\bf y}]=[{\bf x}\;\;{\bf
y}]\;\left[\begin{array}{cc}
\cos(\alpha) & -\sin(\alpha) \\
\sin(\alpha) & \cos(\alpha)\end{array}\right]\,.
\edm

\noindent For $f(\alpha)=\|\tilde{\bf x}\|^2=\tilde{\bf x}^T\tilde{\bf x}$ we have ${\partial f}/{\partial\alpha}=2\,\tilde{\bf x}^T\tilde{\bf y}$.

\end{lem}

\noindent {\bf Proof.} We write $f(\alpha)=\sum_{i=1}^p
(\cos(\alpha)\,x_i+\sin(\alpha)\,y_i)^2$. The first derivative is obtained as 
\bdm
\frac{\partial f}{\partial\alpha}=
2\,\sum_{i=1}^p
(\cos(\alpha)\,x_i+\sin(\alpha)\,y_i)\,(-\sin(\alpha)\,x_i+\cos(\alpha)\,y_i
)= 2\,\sum_{i=1}^p \tilde{x}_i\,\tilde{y}_i = 2\,\tilde{\bf
x}^T\tilde{\bf y}\,. \edm
\ep

\noindent Consider rotating rows $i$ and $j$ of each slice of $\wt{\cal Z}$, with $1\le i\le R_1$ and $R_1+1\le j\le I$. This is done by premultiplying $\wt{\bf S}^T$ by a rotation matrix ${\bf Q}$ that is equal to ${\bf I}_I$ except 
\be
q_{ii}=q_{jj}=\cos(\alpha)\,,\quad\quad\quad
q_{ji}=-q_{ij}=\sin(\alpha)\,.
\ee

\noindent By Lemma~\ref{lem-statorth}, the rotation can increase the norm of the subtensor unless $\tilde{\bf x}$ and $\tilde{\bf y}$ are orthogonal, where $\tilde{\bf x}$ contains the $i$th rows of ${\bf G}_k$, $k=1,\ldots,R_3$, stacked below each other as column vectors, and $\tilde{\bf y}$ analogously contains the $(j-R_1)$th rows of ${\bf N}_k$, $k=1,\ldots,R_3$. This yields (\ref{eq-statM-row}) as the first-order condition for all row rotations together. Analogously, the first-order condition for all column rotations together equals (\ref{eq-statM-col}). For rotations of slices we obtain the first-order condition (\ref{eq-statM-slice}). 

We use the first-order conditions (\ref{eq-statM-row})--(\ref{eq-statM-slice}) in the proof of our main result in section 4. Note that the first-order conditions only imply restrictions on ${\bf S},{\bf T},{\bf U}$ and not on the full matrices $\wt{\bf S},\wt{\bf T},\wt{\bf U}$. 
Indeed, the additional columns in the latter matrices correspond to rotations of rows, columns, or slices outside the $R_1\times R_2\times R_3$ subtensor.

Using results from algebraic geometry, the following lemma is obtained. Its proof is found in the appendix. The notion of ``almost all ${\cal Z}\in\R^{I\times J\times K}$ satisfy property P'' is used in the sense that the Lebesgue measure is zero for the set of ${\cal Z}$ not satisfying P.

\begin{lem}
\label{lem-M-AG}
For almost all ${\cal Z}\in\R^{I\times J\times K}$ problem $(\ref{prob-Mult})$ has a unique minimizer ${\cal X}$ and a finite number of stationary points ${\cal X}_i$, i.e., with corresponding ${\bf S}_i,{\bf T}_i,{\bf U}_i$ satisfying the first-order conditions $(\ref{eq-statM-row})$--$(\ref{eq-statM-slice})$.
\end{lem}

\noindent {\bf Proof.} See Appendix A.
\ep

\section{Finding a best approximation from $\ol{S}_2(I,J,K)$}
\setcounter{equation}{0}
Here we discuss a parameterization of $\ol{S}_2(I,J,K)$ and derive first-order conditions for problem (\ref{prob-CPDc}) with $R=2$. In the sequel we make use of two classifications of tensors in $\R^{2\times 2\times 2}$, both of which can be found in Appendix D. One distinguishes interior, boundary, and exterior points of $S_2(2,2,2)$ and is due to \cite{Ste-arxiv-GSD}. The other is the classification of \cite{DSL} of $\R^{2\times 2\times 2}$ into eight distinct orbits under nonsingular transformations. We start with following lemma.

\begin{lem}
\label{lem-sets}
For $R\le\min(I,J,K)$ we have the following results.
\begin{itemize}
\item[$(i)$] Any ${\cal Y}\in\ol{S}_R(I,J,K)$ can be written as ${\cal Y}=({\bf S},{\bf T},{\bf U})\cdot{\cal H}$, with columnwise orthogonal ${\bf S}\in\R^{I\times R}$, ${\bf T}\in\R^{J\times R}$, and ${\bf U}\in\R^{K\times R}$, and ${\cal H}\in\ol{S}_R(R,R,R)$ having upper triangular slices. Moreover, ${\cal Y}\in S_R(I,J,K)$ if and only if ${\cal H}\in S_R(R,R,R)$.
\item[$(ii)$] $\ol{S}_R(I,J,K)\subseteq M_{(R,R,R)}(I,J,K)$.
\item[$(iii)$] For any ${\cal Y}=({\bf S},{\bf T},{\bf U})\cdot{\cal H}$, with columnwise orthogonal ${\bf S}\in\R^{I\times 2}$, ${\bf T}\in\R^{J\times 2}$, and ${\bf U}\in\R^{K\times 2}$, and ${\cal H}\in\R^{2\times 2\times 2}$ having upper triangular slices, we have ${\cal Y}\in\ol{S}_2(I,J,K)$.
\end{itemize}
\end{lem}

\noindent {\bf Proof.} The proof of $(i)$ can be found in \cite[lemma 3.2]{Ste-limit} and uses \cite[theorem 5.2]{DSL}. Since $({\bf S},{\bf T},{\bf U})\cdot{\cal H}$ in $(i)$ has multilinear rank at most $(R,R,R)$, statement $(ii)$ follows from $(i)$. 

The proof of $(iii)$ is as follows. We have rank$({\cal Y})=$ rank$({\cal H})$ and ${\cal Y}\in\ol{S}_2(I,J,K)$ if and only if ${\cal H}\in\ol{S}_2(2,2,2)$. We use the classification of interior, boundary, and exterior points of $S_2(2,2,2)$ in Proposition~\ref{p-1} in Appendix D.

Let ${\cal H}\in\R^{2\times 2\times 2}$ have upper triangular slices ${\bf H}_1$ and ${\bf H}_2$. If no linear combination of the slices is nonsingular, then ${\cal H}$ is a boundary point of $S_2(2,2,2)$ and ${\cal H}\in\ol{S}_2(2,2,2)$ follows. Next suppose that a linear combination (which is also upper triangular) of the slices exists that is nonsingular. Without loss of generality we assume that ${\bf H}_1$ is nonsingular. If ${\bf H}_2{\bf H}_1^{-1}$ has real eigenvalues, then ${\cal H}\in\ol{S}_2(2,2,2)$. This is true, since ${\bf H}_2{\bf H}_1^{-1}$ is upper triangular. If ${\bf H}_2{\bf H}_1^{-1}$ has complex eigenvalues, then ${\cal H}\notin\ol{S}_2(2,2,2)$. However, ${\bf H}_2{\bf H}_1^{-1}$ is upper triangular and has real eigenvalues. This completes the proof of $(iii)$.
\ep

\noindent From Lemma~\ref{lem-sets} we obtain that
\begin{eqnarray}
\label{eq-S2}
\ol{S}_2(I,J,K)=\{{\cal Y}\in\R^{I\times J\times K}:& \;{\cal Y}=({\bf S},{\bf T},{\bf U})\cdot{\cal H}\,,\;{\rm with}\;{\bf S}^T{\bf S}={\bf T}^T{\bf T}={\bf U}^T{\bf U}={\bf I}_2 \nonumber \\
& {\rm and}\;{\bf H}_k\;{\rm upper\;triangular}\,,k=1,2\}\,.
\end{eqnarray}

\noindent An analogous result is proven in \cite{RoG}, who set
\be
\label{eq-222-RoG}
[{\bf H}_1\;|\;{\bf H}_2]=\left[\begin{array}{cc|cc}
h_{111} & h_{121} & 0 & h_{122} \\
0 & h_{221} & 0 & h_{222}\end{array}\right]\,,
\ee

\noindent which can always be obtained via an orthogonal transformation $({\bf I}_2,{\bf I}_2,{\bf U})\cdot{\cal H}$ with ${\cal H}$ having upper triangular slices. It is shown in \cite{RoG} that ${\cal H}$ in (\ref{eq-222-RoG}) has rank 3 if and only if $h_{222}=0$ (and $h_{111}\neq 0$, $h_{221}\neq 0$, $h_{122}\neq 0$). 

Next, we consider the problem of finding a best approximation of ${\cal Z}\in\R^{I\times J\times K}$ from $\ol{S}_2(I,J,K)$ using the parameterization (\ref{eq-S2}). Analogous to finding a best rank-$(R_1,R_2,R_3)$ approximation, problem (\ref{prob-CPDc}) is equivalent to maximizing the Frobenius norm of the upper triangular parts of the slices of the $2\times 2\times 2$ tensor $({\bf S}^T,{\bf T}^T,{\bf U}^T)\cdot{\cal Z}$ over $({\bf S},{\bf T},{\bf U})$ and setting ${\cal H}$ equal to the upper triangular parts of $({\bf S}^T,{\bf T}^T,{\bf U}^T)\cdot{\cal Z}$. This in turn is equivalent to maximizing the Frobenius norm of the upper triangular parts of the slices of any $2\times 2\times 2$ subtensor of $(\wt{\bf S}^T,\wt{\bf T}^T,\wt{\bf U}^T)\cdot{\cal Z}$ over orthogonal $\wt{\bf S}\in\R^{I\times I}$, $\wt{\bf T}\in\R^{J\times J}$, and $\wt{\bf U}\in\R^{K\times K}$, with ${\cal H}$ equal to the upper triangular parts of the transformed subtensor, and $({\bf S},{\bf T},{\bf U})$ equal to the corresponding columns of $(\wt{\bf S},\wt{\bf T},\wt{\bf U})$. We work with this formulation of problem (\ref{prob-CPDc}) for $R=2$, where we take the subtensor that has indices $(i,j,k)$ with $i=1,2$, $j=1,2$, and $k=1,2$. 

The alternating least squares (ALS) algorithm derived in \cite{Rocc} is used in \cite{RoG} to find a best approximation from $\ol{S}_2(I,J,K)$ using (\ref{eq-S2}) under the restriction (\ref{eq-222-RoG}). Alternatively, (\ref{eq-S2}) can be used without the restriction (\ref{eq-222-RoG}) and an algorithm iterating over Givens rotations can be applied. We derive first-order conditions for the latter problem analogous to (\ref{eq-statM-row})--(\ref{eq-statM-slice}) for a best rank-$(R_1,R_2,R_3)$ approximation. Again, let the slices of $\wt{\cal Z}=(\wt{\bf S}^T,\wt{\bf T}^T,\wt{\bf U}^T)\cdot{\cal Z}$ be partitioned as $\left[\begin{array}{cc}{\bf G}_k & {\bf L}_k\\ {\bf N}_k & {\bf M}_k\end{array}\right]$, where ${\bf G}_k\in\R^{2\times 2}$, ${\bf L}_k\in\R^{2\times (J-2)}$, ${\bf N}_k\in\R^{(I-2)\times 2}$, and ${\bf M}_k\in\R^{(I-2)\times (J-2)}$, $k=1,\ldots,K$. Hence, in an optimal solution the core tensor ${\cal H}$ has slices equal to the upper triangular parts of ${\bf G}_k$, $k=1,2$.

We obtain first-order conditions analogous to (\ref{eq-statM-row})--(\ref{eq-statM-slice}), with $(R_1,R_2,R_3)=(2,2,2)$ and slices ${\bf G}_1$ and ${\bf G}_2$ upper triangular. That is,
\be
\label{eq-statS-row}
\left[\begin{array}{cc|cc} g_{111} & g_{121} & g_{112} & g_{122}\\
0 & g_{221} & 0 & g_{222}\end{array}\right]\;[{\bf N}_1\;{\bf N}_2]^T={\bf O}\,,
\ee
\be
\label{eq-statS-col}
\left[\begin{array}{cc} g_{111} & g_{121} \\ 0 & g_{221} \\[1mm] \hline 
g_{112} & g_{122}\\ 0 & g_{222}\end{array}\right]^T\;
\left[\begin{array}{c}{\bf L}_1\\ {\bf L}_2\end{array}\right]={\bf O}\,,
\ee
\be
\label{eq-statS-slice}
\left[\begin{array}{c|c} g_{111} & g_{112} \\ 0 & 0 \\ g_{121} & g_{122}\\
g_{221} & g_{222}\end{array}\right]^T\;
\left[\begin{array}{c|c|c} g_{113} & \cdots & g_{11K} \\ 0 & \cdots & 0 \\
g_{123} & \cdots & g_{12K} \\ g_{223} & \cdots & g_{22K}\end{array}\right]={\bf O}\,.
\ee

\noindent Additionally, using Lemma~\ref{lem-statorth}, we need to consider rotations of rows 1 and 2 and columns 1 and 2, since these too transfer mass to the upper triangular parts of ${\bf G}_1$ and ${\bf G}_2$. For rows 1 and 2 we obtain the condition
\be
\label{eq-statS-row12}
\left(\begin{array}{c} g_{211}\\ g_{212}\end{array}\right)^T
\left(\begin{array}{c} g_{111}\\ g_{112}\end{array}\right)=0\,.
\ee

\noindent For columns 1 and 2 we obtain the condition
\be
\label{eq-statS-col12}
\left(\begin{array}{c} g_{211}\\ g_{212}\end{array}\right)^T
\left(\begin{array}{c} g_{221}\\ g_{222}\end{array}\right)=0\,.
\ee

\noindent We use the first-order conditions (\ref{eq-statS-row})--(\ref{eq-statS-col12}) in the proof of our main result in section 4. Note that the first-order conditions only imply restrictions on ${\bf S},{\bf T},{\bf U}$ (with ${\cal G}=({\bf S}^T,{\bf T}^T,{\bf U}^T)\cdot{\cal Z}$) and not on the full matrices $\wt{\bf S},\wt{\bf T},\wt{\bf U}$. Analogous to Lemma~\ref{lem-M-AG} we have the following result. 

\begin{lem}
\label{lem-S-AG}
For almost all ${\cal Z}\in\R^{I\times J\times K}$ problem $(\ref{prob-CPDc})$ has a unique minimizer ${\cal X}$ and a finite number of stationary points ${\cal X}_i$, i.e., with corresponding ${\bf S}_i,{\bf T}_i,{\bf U}_i$ satisfying the first-order conditions $(\ref{eq-statS-row})$--$(\ref{eq-statS-col12})$.
\end{lem}

\noindent {\bf Proof.} See Appendix A.
\ep

\noindent For later use, we state the following results.

\begin{lem}
\label{lem-S-rank222}
For almost all ${\cal Z}\in\R^{I\times J\times K}$ with ${\cal Z}\notin\ol{S}_2(I,J,K)$ any local minimizer ${\cal X}$ of problem $(\ref{prob-CPDc})$ has {\rm mrank}$({\cal X})=(2,2,2)$.
\end{lem}

\noindent {\bf Proof.} See Appendix B.
\ep

\begin{lem}
\label{lem-optrank}
Let the set $W_R(I,J,K)$ denote either $S_R(I,J,K)$, $\ol{S}_R(I,J,K)$, or $M_{(R,R,R)}(I,J,K)$, and let ${\cal X}$ be a best approximation of ${\cal Z}\notin W_R(I,J,K)$ from the set $W_R(I,J,K)$. Then {\rm rank}$({\cal X})\ge R$.
\end{lem}

\noindent {\bf Proof.} The proof for $S_R(I,J,K)$ is easy: if rank$({\cal X})<R$, then a rank-1 term can be added to ${\cal X}$ to obtain a better approximation of ${\cal Z}$ \cite[lemma 8.2]{DSL}. Since $S_R(I,J,K)\subseteq\ol{S}_R(I,J,K)\subseteq M_{(R,R,R)}(I,J,K)$ (see Lemma~\ref{lem-sets}), the same argument can be used for $\ol{S}_R(I,J,K)$ and $M_{(R,R,R)}(I,J,K)$.
\ep

\section{Existence of a best rank-2 approximation}
\setcounter{equation}{0}
Here we relate existence of a best rank-2 approximation to the rank of local minima of the best rank-(2,2,2) approximation problem. Since $S_2(I,J,K)\subseteq M_{(2,2,2)}(I,J,K)$ (Lemma~\ref{lem-sets}), it follows that if a best rank-(2,2,2) approximation has rank 2, then it is also a best rank-2 approximation. Our main result is a partial converse of this. Recall that a best rank-(2,2,2) approximation has rank equal to the rank of its core tensor ${\cal G}\in\R^{2\times 2\times 2}$, which has rank 2 or 3 (Lemma~\ref{lem-optrank}, and the fact that the maximal rank equals 3; see Appendix D). We prove the following.
\begin{thm}
\label{t-main}
Let ${\cal Z}\in\R^{I\times J\times K}$ have rank larger than $2$. Let ${\cal X}$ be a locally best approximation of ${\cal Z}$ from $\ol{S}_2(I,J,K)$, with {\rm rank}$({\cal X})=2$. The following statements hold.
\begin{itemize}
\item[$(i)$] ${\cal X}$ is a stationary point of the best {\rm rank}-$(2,2,2)$ approximation problem for ${\cal Z}$.
\item[$(ii)$] If ${\cal X}$ is an interior point of $S_2(I,J,K)$, then it is a local minimizer in the best {\rm rank}-$(2,2,2)$ approximation problem for ${\cal Z}$.
\end{itemize}
\end{thm}

\noindent {\bf Proof.} First, we prove $(i)$. We consider a locally best approximation ${\cal X}$ of ${\cal Z}$ from the set $\ol{S}_2(I,J,K)$ parameterized as in (\ref{eq-S2}). We write ${\cal X}=({\bf S},{\bf T},{\bf U})\cdot{\cal H}$ with columnwise orthogonal ${\bf S}\in\R^{I\times 2}$, ${\bf T}\in\R^{J\times 2}$, and ${\bf U}\in\R^{K\times 2}$, and ${\cal H}\in\R^{2\times 2\times 2}$ having upper triangular slices. We assume that rank$({\cal X})=$ rank$({\cal H})=2$ and consider all possibilities for ${\cal H}$ in Proposition~\ref{p-1}.

Suppose first that slice ${\bf H}_1$ is nonsingular (or, equivalently, that a linear combination of the slices of ${\cal H}$ is nonsingular). Proposition~\ref{p-1} implies that ${\bf H}_2{\bf H}_1^{-1}$ has real eigenvalues and two linearly independent eigenvectors. Recall that ${\bf H}_k$ is equal to the upper triangular part of ${\bf G}_k$ in the first-order conditions (\ref{eq-statS-row})--(\ref{eq-statS-col12}). We have
\be
\label{eq-H12}
{\bf H}_2{\bf H}_1^{-1}=\left[\begin{array}{cc}
g_{111}^{-1}\,g_{112} & g_{221}^{-1}\,(g_{122}-g_{111}^{-1}\,g_{121}\,g_{112}) \\
0 & g_{221}^{-1}\,g_{222}\end{array}\right]\,,
\ee

\noindent where $g_{111}\neq 0$ and $g_{221}\neq 0$ since ${\bf H}_1$ is nonsingular.
Let the eigenvalues of ${\bf H}_2{\bf H}_1^{-1}$ be distinct. By (\ref{eq-H12}) this is equivalent to the vectors $(g_{111}\;g_{112})$ and $(g_{221}\;g_{222})$ being linearly independent. From (\ref{eq-statS-row12})--(\ref{eq-statS-col12}) we then obtain $g_{211}=g_{212}=0$. Hence, ${\bf G}_k={\bf H}_k$ are upper triangular for $k=1,2$, and first-order conditions (\ref{eq-statS-row})--(\ref{eq-statS-slice}) are identical to first-order conditions (\ref{eq-statM-row})--(\ref{eq-statM-slice}). 

Suppose next that ${\bf H}_1$ is nonsingular and ${\bf H}_2{\bf H}_1^{-1}$ has two identical real eigenvalues with two linearly independent eigenvectors. From (\ref{eq-H12}) it follows that the vectors $(g_{111}\;g_{112})$ and $(g_{221}\;g_{222})$ are proportional and that $g_{111}\,g_{122}=g_{121}\,g_{112}$. The latter implies that the vectors $(g_{111}\;g_{112})$ and $(g_{121}\;g_{122})$ are also proportional. (The first-order conditions (\ref{eq-statS-row12})--(\ref{eq-statS-col12}) imply that the vector $(g_{211}\;g_{212})$ is orthogonal to the three proportional vectors mentioned.)
It is now possible to rotate rows 1 and 2 or columns 1 and 2 such that ${\bf G}_k=\left[\begin{array}{cc} *&0\\ *&*\end{array}\right]$, $k=1,2$. Swapping rows 1 and 2, followed by swapping columns 1 and 2, then yields upper triangular ${\bf G}_k$, $k=1,2$. Hence, we obtain a better approximation from $\ol{S}_2(I,J,K)$ unless ${\bf G}_k$, $k=1,2$, were already upper triangular before the rotations. As above we obtain 
${\bf G}_k={\bf H}_k$, $k=1,2$.

Finally, suppose that no linear combination of ${\bf H}_1$ and ${\bf H}_2$ exists that is nonsingular. Since ${\bf H}_1$ and ${\bf H}_2$ are upper triangular, it follows that $g_{111}=g_{112}=0$ or $g_{221}=g_{222}=0$ or both. If $g_{111}=g_{112}=0$, then a rotation of rows 1 and 2 can move nonzero mass to the (1,1) entries unless $g_{211}=g_{212}=0$. Likewise, if $g_{221}=g_{222}=0$, then a rotation of columns 1 and 2 can move nonzero mass to the (2,2) entries unless $g_{221}=g_{222}=0$. Hence, in both cases ${\bf G}_k={\bf H}_k$ are upper triangular for $k=1,2$. Hence, in all possible cases ${\cal X}$ is a stationary point of the best {\rm rank}-$(2,2,2)$ approximation. This completes the proof of $(i)$.

Next, we prove $(ii)$. Local minimizer ${\cal X}=({\bf S},{\bf T},{\bf U})\cdot{\cal H}$ is an interior point of $S_2(I,J,K)$ if and only if ${\cal H}$ is an interior point of $S_2(2,2,2)$. Note that mrank$({\cal X})=$ mrank$({\cal H})=(2,2,2)$ according to Proposition~\ref{p-1}. We use the fact that generic tensors in $\R^{2\times 2\times 2}$ have rank 2 or 3, both on sets of positive Lebesgue measure (Appendix D), to obtain that for ${\cal Y}$ in a small neighborhood of ${\cal X}$ in $M_{(2,2,2)}(I,J,K)$ we have rank$({\cal Y})=2$. It follows that $\|{\cal Z}-{\cal Y}||_F^2\ge\|{\cal Z}-{\cal X}\|^2_F$, which implies that ${\cal X}$ is a local minimizer of the best rank-$(2,2,2)$ approximation problem. This completes the proof.
\ep

\noindent We have the following corollary to Theorem~\ref{t-main}.

\begin{cor}
\label{cor-main}
For almost all ${\cal Z}\in\R^{I\times J\times K}$ with rank$({\cal Z})>2$, if all locally best approximations ${\cal X}_i$ of ${\cal Z}$ from $M_{(2,2,2)}(I,J,K)$ satisfy {\rm rank}$({\cal X}_i)=3$, then ${\cal Z}$ does not have a best rank-$2$ approximation.
\end{cor}

\noindent {\bf Proof.} Recall from Lemma~\ref{lem-M-AG} that there are finitely many locally best approximations ${\cal X}_i$ of ${\cal Z}$ from $M_{(2,2,2)}(I,J,K)$. Theorem~\ref{t-main} $(ii)$, together with the conditions of the corollary, implies that ${\cal Z}$ has no best rank-2 approximation that is an interior point of $S_2(I,J,K)$. It remains to consider the possibility that a best rank-2 approximation ${\cal X}$ of ${\cal Z}$ is a boundary point of $S_2(I,J,K)$. From Lemma~\ref{lem-optrank} it follows that rank$({\cal X})=2$. From the classification of tensors in $\R^{2\times 2\times 2}$ into eight orbits (Appendix D) it follows that ${\cal X}$ has multilinear rank (1,2,2), (2,1,2), or (2,2,1). This possibility is excluded by Lemma~\ref{lem-S-rank222}.
\ep

\noindent We now have the following options for determining whether ${\cal Z}\in\R^{I\times J\times K}$ has a best rank-2 approximation or not, and to obtain a best rank-2 approximation when it exists. 
\begin{itemize}
\item[(A)] Compute the best approximation ${\cal X}=({\bf S},{\bf T},{\bf U})\cdot{\cal H}$ from the set $\ol{S}_2(I,J,K)$ using the parameterization (\ref{eq-S2}) under the restriction (\ref{eq-222-RoG}), and the algorithm of \cite{Rocc}. As shown in \cite{RoG} the solution has rank 3 if and only if $h_{222}=0$ (and $h_{111}\neq 0$, $h_{221}\neq 0$, $h_{122}\neq 0$).
\item[(B)] Compute the best approximation ${\cal X}=({\bf S},{\bf T},{\bf U})\cdot{\cal H}$ from the set $\ol{S}_2(I,J,K)$ using the parameterization (\ref{eq-S2}) by iterating over Givens rotations. The rank of ${\cal H}$ follows from the criteria in Proposition~\ref{p-1}.
\item[(C)] Compute all local minima ${\cal X}=({\bf S},{\bf T},{\bf U})\cdot{\cal G}$ of the best rank-(2,2,2) approximation problem. The rank of ${\cal G}$ for a best rank-(2,2,2) approximation determines whether a best rank-2 approximation has been found, and all ranks of ${\cal G}$ may determine whether it exists (Theorem~\ref{t-main}).
\end{itemize}

\noindent For (A) we need to check whether $h_{222}$ equals zero or not, and for (B) we need to verify whether the eigenvalues of ${\bf H}_2{\bf H}_1^{-1}$ are identical or not. In practice, these criteria can be verified numerically only be setting some tolerance. For large tensors the tolerance may need to be larger as well, and the stopping criterion of the algorithm to find a best approximation from $\ol{S}_2(I,J,K)$ may need to be more strict. For (C) we need to verify whether ${\bf G}_2{\bf G}_1^{-1}$ has real or complex eigenvalues (Proposition~\ref{p-1}), which is numerically more reliable. Note, however, that (C) does not yield a conclusion regarding existence of a best rank-2 approximation when local rank-(2,2,2) minima are encountered of both rank 2 and rank 3, but the best rank-(2,2,2) approximation has rank 3.\\

\noindent When using (C) we do not expect to encounter any locally best rank-(2,2,2) approximations that are a boundary points of $\ol{S}_2(I,J,K)$ (i.e., with ${\bf G}_2{\bf G}_1^{-1}$ having identical real eigenvalues, see Proposition~\ref{p-1}). This is formally proven in the lemma below.

\begin{lem}
\label{lem-M-interior}
For almost all ${\cal Z}\in\R^{I\times J\times K}$ with {\rm mrank}$({\cal Z})>(2,2,2)$ any locally best {\rm rank}-$(2,2,2)$ approximation ${\cal X}$ has {\rm mrank}$({\cal X})=(2,2,2)$ and is not a boundary point of $\ol{S}_2(I,J,K)$.
\end{lem}

\noindent {\bf Proof.} See Appendix C.
\ep

\newpage
\section{Simulations}
\setcounter{equation}{0}
We conduct a simulation study to compare the criteria (A) and (C) for existence of a best rank-2 approximation for randomly sampled tensors. For (A) we use the ALS algorithm of \cite{Rocc} as is also done in \cite{RoG}. For (C) we use the ALS algorithm of \cite{KrL} as described in \cite{LMV00}. For both algorithms we use convergence criterion $10^{-9}$ for the relative decrease of $\|{\cal Z}-{\cal Y}\|_F^2$. The entries of tensor ${\cal Z}$ are sampled independently from the standard normal distribution. For each ${\cal Z}$ we run both algorithms 40 times: 39 times with random starting values and 1 time with starting values computed from the SVDs of the three matrix unfoldings of ${\cal Z}$ as suggested in \cite{KrL}. Hence, for each ${\cal Z}$ and each algorithm we obtain 40 local minima. We determine the rank of these minima via checking whether $h_{222}=0$ for (A) and checking the eigenvalues of ${\bf G}_2{\bf G}_1^{-1}$ for (C). All local minima for (A) satisfy $h_{111}\neq 0$, $h_{221}\neq 0$, and $h_{122}\neq 0$, which implies they have multilinear rank (2,2,2). Hence, no counterexamples to Lemma~\ref{lem-S-rank222} were found. All local minima for (C) have nonsingular ${\bf G}_1$ and the eigenvalues of ${\bf G}_2{\bf G}_1^{-1}$ are distinct, which is in line with Lemma~\ref{lem-M-interior}.

To illustrate the difference between small and large tensors, we consider two sizes of tensors: $4\times 4\times 4$ and $30\times 10\times 5$. The latter size was also used in the simulation study in \cite{RoG}. For each size we generate 1000 random tensors ${\cal Z}$. The ALS algorithms of \cite{Rocc} and \cite{KrL} yield a monotonically decreasing sequence $\|{\cal Z}-{\cal Y}^{(n)}\|^2_F$. To check whether they terminate in local minima (and not saddle points) we compute the eigenvalues of the Hessian matrix of second-order derivatives corresponding to the maximization of $\|({\bf S}^T,{\bf T}^T,{\bf U}^T)\cdot{\cal Z}\|^2_F$ (or the upper triangular parts thereof for $\ol{S}_2(I,J,K)$). Here, the variables of the problem are the rotation angles of the Givens rotations parameterizing orthogonal ${\bf Q}_1\in\R^{I\times I}$, ${\bf Q}_2\in\R^{J\times J}$, and ${\bf Q}_3\in\R^{K\times K}$, and the objective function is taken as $\|({\bf S}^T{\bf Q}_1,{\bf T}^T{\bf Q}_2,{\bf U}^T{\bf Q}_3)\cdot{\cal Z}\|^2_F$ (or the upper triangular parts thereof for $\ol{S}_2(I,J,K)$). The Hessian is computed symbolically and then evaluated at ${\bf Q}_1={\bf I}_I$, ${\bf Q}_2={\bf I}_J$, and ${\bf Q}_3={\bf I}_K$. The largest eigenvalue of the Hessian should be close to zero for a local maximum. Note that the Hessian always has zero eigenvalues, since rotations of rows 1 and 2, columns 1 and 2, and slices 1 and 2 do not affect the objective function of finding a best rank-(2,2,2) approximation. For finding a best approximation from $\ol{S}_2(I,J,K)$, (more complicated) combinations of these rotations can be found that do not affect the objective function. Since the Hessian is computed symbolically, which takes a lot of time, we compute the Hessian only for 100 random $4\times 4\times 4$ tensors and one run with random starting values per tensor of both ALS algorithms. The largest eigenvalue of the Hessian for all 100 runs is $0.00033$ for the algorithm of \cite{Rocc} and $0.00025$ for the algorithm of \cite{KrL}. Hence, we have found only local minima as solutions produced by the algorithms in these runs.

Next, we report the number of distinct local minima found by each algorithm for the same tensor ${\cal Z}$ (for all 1000 generated tensors). We consider two solutions ${\cal X}_1$ and ${\cal X}_2$ as distinct when $\|{\cal X}_1-{\cal X}_2\|_F^2>\epsilon$ for some small threshold $\epsilon>0$. For $4\times 4\times 4$ tensors we take $\epsilon=0.001$ and for $30\times 10\times 5$ we take $\epsilon=0.1$. Table~\ref{tab-minfreq} displays the numbers of distinct local minima found for both algorithms. For $4\times 4\times 4$ tensors up to 6 local minima are found, and for $30\times 10\times 5$ tensors up to 22 local minima are found. These numbers do not change much for slightly different thresholds $\epsilon$. \\

\begin{table}[h]
\begin{center}
\begin{tabular}{|c|cccccccccc|r|}
\hline & \multicolumn{10}{c|}{number of distinct local minima} & \\
\hline 
problem & 1 & 2 & 3 & 4 & 5 & 6 & 7 & 8 & 9 & 10--22 &  total \\
\hline
$\ol{S}_2(4,4,4)$ & 301 & 419 & 213 & 55 & 10 &  2 &  0 & 0 & 0&0&2060 \\
$M_{(2,2,2)}(4,4,4)$ & 379 & 418 & 171 & 30 &  1  & 1 & 0 & 0 & 0&0&1859\\
\hline
$\ol{S}_2(30,10,5)$ & 0  & 9 & 36 & 66 & 91 & 121 & 138 & 145 & 118   & 276  & 8019 \\
$M_{(2,2,2)}(30,10,5)$ & 2 & 11 & 40 & 81 & 104 & 139 & 160 & 147 & 97  & 219 & 7567 \\
\hline 
\end{tabular}
\end{center}
\caption{Number of tensors ${\cal Z}$ with number of distinct local minima as indicated per column. The last column contains the total number of distinct local minima found.} 
\label{tab-minfreq}
\end{table}

\noindent Depending on the eigenvalues of ${\bf G}_2{\bf G}_1^{-1}$ for the local minima of the best rank-(2,2,2) problem, we partition the 1000 tensors into four subsets: all local minima having real eigenvalues (rank 2, referred to as `all real'), all local minima having complex eigenvalues (rank 3, referred to as `all complex'), both real and complex eigenvalues occur but the best minimum has real eigenvalues (`mixed real'), and both real and complex eigenvalues occur but the best minimum has complex eigenvalues (`mixed complex'). In Table~\ref{tab-minima} we report for each subset the number of rank-(2,2,2) minima with real eigenvalues that are also found as minima in the approximation from $\ol{S}_2(I,J,K)$, and their associated values for $|h_{222}|$. \\

\begin{table}[h]
\begin{center}
\begin{tabular}{|c||c|c|c|c|c|c|}
\hline 
eig$({\bf G}_2{\bf G}_1^{-1})$ &  \# tensors & \# real & \# shared & \# shared  & 
$|h_{222}|$ in & $|h_{222}|$ in \\[-2mm]
&  & min & min & best min & shared min & best min \\
\hline
all real & 644 & 1167 & 1155 & 644 & 0.25--5.28 & 0.25--5.22\\
mixed real & 114 & 169 & 168 & 114 & 0.51--5.05 & 1.09--5.05 \\
all complex & 121 & - & - & - & - & 0.00--0.16 \\
mixed complex & 121 &  145 &  143 & - & 0.86--5.36 & 0.00--4.12 \\
\hline
all real & 209 & 1309 & 1178 & 208 & 0.05--10.36 & 0.59--10.30 \\
mixed real & 542 & 3554 & 3101 & 542 & 0.29--10.42 & 0.64--10.04 \\
all complex & 1 &- &- &- &- & 0.03 \\
mixed complex & 248 & 1410 & 1159 & - & 0.28--10.33 & 0.00--9.47 \\
\hline
\end{tabular}
\end{center}
\caption{For each subset of tensors are given: the number of tensors ${\cal Z}$, number of rank-2 minima of the best rank-(2,2,2) problem, number of shared such rank-2 minima, number of best such minima that are shared, and minimal and maximal $|h_{222}|$ values for shared and best minima of the approximation from $\ol{S}_2(I,J,K)$. Top rows: $4\times 4\times 4$ tensors. Bottom rows: $30\times 10\times 5$ tensors.} 
\label{tab-minima}
\end{table}

\begin{figure}[t]
\begin{center}
\includegraphics[width=8cm,height=6cm]{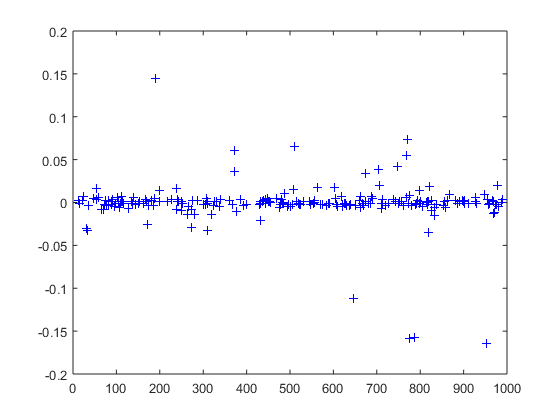}
\includegraphics[width=8cm,height=6cm]{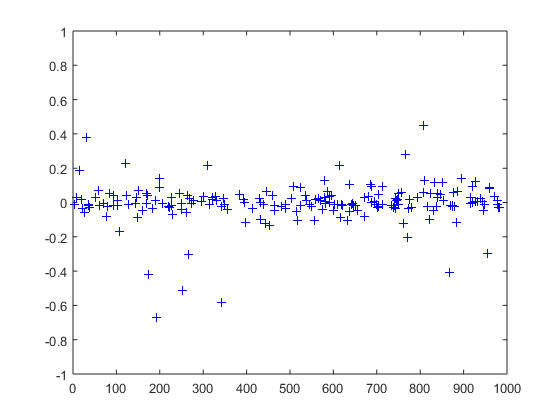}
\caption{Small values of $h_{222}$ for the best approximation from $\ol{S}_2(I,J,K)$ for all 1000 tensors. Left: $4\times 4\times 4$ tensors. Right: $30\times 10\times 5$ tensors.} 
\label{fig-h222} 
\end{center} 
\end{figure}

\noindent For the $4\times 4\times 4$ tensors all four subsets are nonempty and 99 percent of the rank-(2,2,2) minima with real eigenvalues are also found as minima of the  $\ol{S}_2(4,4,4)$ problem. Moreover, all best of such rank-(2,2,2) minima are also found as best minima in the $\ol{S}_2(4,4,4)$ problem. The values of $|h_{222}|$ for the shared minima should be nonzero and their minimal value is 0.25. For the `all real' and `mixed real' subsets both algorithms find the same best approximation, which is in line with theory. For the `all complex' subset Corollary~\ref{cor-main} implies that no best rank-2 approximation exists. The values of $|h_{222}|$ of the best minima in the $\ol{S}_2(4,4,4)$ problem are indeed small, with a maximum of 0.16. However, there is not much difference with the smallest value of 0.25 for the rank-2 minima. For the `mixed complex' subset a best rank-2 approximation may or may not exist and the $|h_{222}|$ values of the best minima in the $\ol{S}_2(4,4,4)$ problem range from nearly zero to 4.12. To determine whether these minima have rank 2 or rank 3 a threshold with respect to $|h_{222}|$ needs to be specified. In Figure~\ref{fig-h222} the small values of $h_{222}$ for the best minimum from $\ol{S}_2(4,4,4)$ are plotted. As can be seen, determining a suitable threshold may be difficult. Checking for real or complex eigenvalues of ${\bf G}_2{\bf G}_1^{-1}$ is numerically more reliable. 

Next, we discuss the results for $30\times 10\times 5$ tensors. In this case the `all complex' subset consists of one tensor only. Hence, the result of Corollary~\ref{cor-main} does not add much to the analysis of the simulation results. Of the rank-(2,2,2) minima with real eigenvalues more than 86 percent are also found as minima of the $\ol{S}_2(30,10,5)$ problem. Almost all best such minima are also found as best minima in the $\ol{S}_2(30,10,5)$ problem. However, one best rank-(2,2,2) minima is better than the best found approximation from $\ol{S}_2(30,10,5)$, and there is also one best approximation from $\ol{S}_2(30,10,5)$ that is better than the best rank-(2,2,2) approximation. Apparently, for thess ${\cal Z}$ more runs of the algorithms are needed. The values of $|h_{222}|$ for the shared minima should be nonzero, but their minimal value is 0.05. This is an outlier, however, with other values being at least 0.58 for the `all real' subset. In Figure~\ref{fig-h222} it can be seen that for $30\times 10\times 5$ tensors the small values of $h_{222}$ show more variation than for $4\times 4\times 4$ tensors, which makes it harder to choose a threshold for $|h_{222}|$. \\

\newpage
\noindent The simulation results clearly show the added value of computing local rank-(2,2,2) minima to determine whether a best rank-2 approximation exists, and to obtain it when it does. Determining the rank of the $2\times 2\times 2$ core tensor ${\cal G}$ via observing real or complex eigenvalues of $({\bf G}_2{\bf G}_1^{-1})$ is numerically more reliable than deciding whether an obtained value of $h_{222}$ is zero or not. Note that both approaches require running the algorithm a large number of times, either to obtain all local minima, or to make sure the best minimum is obtained. However, the rank-(2,2,2) approach is inconclusive for tensors in the `mixed complex' subset, the size of which is 12 and 25 percent in our simulation study. Hence, for these tensors the only option is to compute a best approximation from $\ol{S}_2(I,J,K)$ and use either method (A) or (B) in section 4.

\section{Results for complex tensors}
\setcounter{equation}{0}
As in the real case, for complex tensors a best rank-$R$ approximation may not exist for $R\ge 2$. The contrived examples in \cite{Bini-1,Paa2,DSL} for real tensors are also valid in the complex case. However, whether nonexistence of a best rank-$R$ approximation holds on a set of positive volume is still an open problem. A related result by \cite{VN} concerns tensors that allow a Schmidt-Eckart-Young (SEY) decomposition: a decomposition into a sum of rank-1 tensors with a best rank-$R$ approximation given by the sum of $R$ of these rank-1 tensors. For matrices the SEY decomposition is given by the SVD \cite{EY}. The result of \cite{VN} is that the set of complex order-$N$ tensors that do not admit an SEY decomposition has positive volume for $N\ge 3$. Since the set of complex tensors that do not have a best rank-$R$ approximation is a subset of the set of complex tensors that do not admit an SEY decomposition, this does not imply that the former set also has positive volume. 

For ${\cal Z}\in\C^{I\times J\times 2}$ the situation is better than in the real case. Namely, it is shown in \cite{LW} that $M_{(R,R,2)}(I,J,2)=\ol{S}_R(I,J,2)$ for $R\le\min(I,J)$, and that
$M_{(R,J,2)}(I,J,2)=\ol{S}_R(I,J,2)$ for $J\le R\le I$. Hence, in these cases computing the best rank-$(R,R,2)$ or rank-$(R,J,2)$ approximation suffices to determine whether a best rank-$R$ approximation exists, and to obtain it when it does. Note that Lemma~\ref{lem-M-AG} also holds for complex tensors (see Appendix A).

For complex $2\times 2\times 2$ tensors the generic rank is 2 and the maximal rank is 3 (see Appendix D). Since any rank-3 tensor in $\C^{2\times 2\times 2}$ can be approximated arbitrarily well by rank-2 tensors, it follows that $\C^{2\times 2\times 2}=\ol{S}_2(2,2,2)$. For ${\cal Y}\in M_{(2,2,2)}(I,J,K)$ we have ${\cal Y}=({\bf S},{\bf T},{\bf U})\cdot{\cal G}$, with ${\cal G}\in\C^{2\times 2\times 2}$. Hence, ${\cal Y}\in\ol{S}_2(I,J,K)$. Since $\ol{S}_2(I,J,K)\subseteq M_{(2,2,2)}(I,J,K)$, we obtain $\ol{S}_2(I,J,K)=M_{(2,2,2)}(I,J,K)$.  We thus obtain the following result.

\begin{lem}
\label{lem-C222}
${\cal Z}\in\C^{I\times J\times K}$ has a best $($complex$)$ {\rm rank-2} approximation if and only if it has a best $($complex$)$ {\rm rank-(2,2,2)} approximation of rank $2$.
\ep
\end{lem}

\noindent In algebraic geometry the set $\ol{S}_2(I,J,K)$ is also known as the second secant variety to the Segre variety and the equality $\ol{S}_2(I,J,K)=M_{(2,2,2)}(I,J,K)$ is shown in \cite{Rai}, also for arbitrary order $N$ of the tensors. In \cite{BB} an algorithm is given to check whether a complex order-$N$ tensor has border rank 2, which makes it a boundary point of $\ol{S}_2(I,J,K)$. The description above and in Appendix D for $N=3$ is simpler, however. Note that our results for real tensors do not follow immediately from the mentioned results for complex tensors.

\section{Discussion}
\setcounter{equation}{0}
In this paper we have proposed to compute locally best rank-(2,2,2) approximations for real order-3 tensors to determine whether a best rank-2 approximation exists, and to obtain it when it does. When a best rank-(2,2,2) approximation has rank 2, it is also a best rank-2 approximation. For generic tensors, we proved that when all locally best rank-(2,2,2) approximations have rank 3, then no best rank-2 approximation exists. Verifying the rank of locally best rank-(2,2,2) approximations boils down to checking whether a $2\times 2$ matrix has real or complex eigenvalues, which is numerically more reliable than checking whether a core entry is zero or not \cite{RoG}. This was demonstrated clearly in our simulation study. However, a drawback of our method is that it yields no conclusion regarding the existence of a best rank-2 approximation when the best rank-(2,2,2) approximation has rank 3, but also local minima of rank 2 are encountered. In that case one still has to compute a best approximation from the closure of the rank-2 set and decide whether it has rank 2 or 3. Although using both algorithms for the same tensor may help to choose a suitable tolerance for $|h_{222}|$ by inspecting the values of $h_{222}$ for local minima found by both algorithms. 

One may wonder whether our approach can be extended to best rank-3 approximations. Since $\ol{S}_3(I,J,K)\subset M_{(3,3,3)}(I,J,K)$ it is still true that if the best rank-(3,3,3) has rank 3, then it is also a best rank-3 approximation. However, since real $3\times 3\times 3$ tensors have generic rank 5 \cite{Com09}, we do not expect a best rank-(3,3,3) approximation to have rank 3 for almost all ${\cal Z}$. Indeed, we were unable to find a counterexample for randomly sampled tensors. To obtain an analogue of Theorem~\ref{t-main} for $\ol{S}_3(I,J,K)$ and $M_{(3,3,3)}(I,J,K)$, we need a parameterization of $\ol{S}_3(I,J,K)$ analogous to (\ref{eq-S2}). From Lemma~\ref{lem-sets} $(i)$ we obtain that ${\cal Y}\in\ol{S}_3(I,J,K)$ can be written as $({\bf S},{\bf T},{\bf U})\cdot{\cal H}$, with ${\cal H}\in\R^{3\times 3\times 3}$ having upper triangular slices. However, to obtain rank$({\cal H})=3$ we also need ${\bf H}_2{\bf H}_1^{-1}$ and ${\bf H}_3{\bf H}_1^{-1}$ to have identical eigenvectors \cite{Ste-limit} (assuming ${\bf H}_1$ is nonsingular). Computing a best approximation from $\ol{S}_3(I,J,K)$ does not seem as simple as finding orthogonal $\wt{\bf S},\wt{\bf T},\wt{\bf U}$ that maximize the Frobenius norm of some part of a subtensor of $(\wt{\bf S}^T,\wt{\bf T}^T,\wt{\bf U}^T)\cdot{\cal Z}$. This makes the link with locally best approximations from $M_{(3,3,3)}(I,J,K)$ more complicated than in the rank-2 case.

Note that we may not replace $M_{(3,3,3)}(I,J,K)$ in the above by either $M_{(3,2,2)}(I,J,K)$ or $M_{(3,3,2)}(I,J,K)$, since both these sets do not contain all rank-3 tensors when $\min(I,J,K)\ge 3$. Indeed, an example is the tensor ${\cal Y}=\sum_{r=1}^3({\bf a}_r\circ{\bf b}_r\circ{\bf c}_r)$ with linearly independent ${\bf a}_r$, $r=1,2,3$, and ${\bf b}_r$, $r=1,2,3$, and ${\bf c}_r$, $r=1,2,3$, which has rank 3 and mrank (3,3,3). If $\min(I,J,K)=2$, then $\ol{S}_3(I,J,K)$ is equal to the GSD set and computing a best GSD approximation yields a best rank-3 approximation when it exists \cite{SDL,Ste-arxiv-GSD}.

\section*{Appendix A: proofs of Lemma~\ref{lem-M-AG} and Lemma~\ref{lem-S-AG}}
\refstepcounter{section}
\setcounter{equation}{0}
\renewcommand{\thesection}{A}
\renewcommand{\theequation}{A.\arabic{equation}}
Lemma~\ref{lem-M-AG} and Lemma~\ref{lem-S-AG} are special cases of results obtained in algebraic geometry. Next, we state their proofs. The set $M_{(R_1,R_2,R_3)}(I,J,K)$ is defined by conditions on the ranks of the three matrix unfoldings of an $I\times J\times K$ tensor. A matrix has rank at most $R$ when all minors of order $R$ (i.e., determinants of $R\times R$ submatrices) are zero. These are polynomial equations in the entries of the matrix. Hence, $M_{(R_1,R_2,R_3)}(I,J,K)$ can be defined in terms of polynomial equations in the entries of an $I\times J\times K$ tensor. It follows that $M_{(R_1,R_2,R_3)}(I,J,K)$ is a variety. The set $\overline{S}_R(I,J,K)$ can be defined as the closure of the image of a polynomial map \cite{Fri12} and, hence, it is a (irreducible) variety. 

To show that almost all ${\cal Z}$ have a unique best approximation from $M_{(R_1,R_2,R_3)}(I,J,K)$ and from $\overline{S}_R(I,J,K)$, we apply \cite[theorem 3.7]{FriSta15}. Requirements in \cite{FriSta15} are that the set is closed and semi-algebraic and the norm is semi-algebraic and differentiable. Since $M_{(R_1,R_2,R_3)}(I,J,K)$ and $\overline{S}_R(I,J,K)$ are closed, a variety is by definition semi-algebraic, and the Frobenius norm is semi-algebraic and differentiable \cite{FriSta15}, this completes the proof of statements $(i)$ of Lemma~\ref{lem-M-AG} and Lemma~\ref{lem-S-AG}. 

Next, we prove that the number of stationary points is finite in the approximation problems from $M_{(R_1,R_2,R_3)}(I,J,K)$ and $\overline{S}_2(I,J,K)$. We begin with $M_{(R_1,R_2,R_3)}(I,J,K)$. It is known that the set of $I\times J$ matrices with rank exactly $R$ is a manifold, i.e., in a neighborhood of each point a homeomorphism to Euclidian space of dimension $n$ exists. Such points are called smooth points and, hence, a manifold consists of only smooth points. Analogous to the matrix result, the set of $I\times J\times K$ tensors with multilinear rank exactly $(R_1,R_2,R_3)$ is also a manifold. After fixing orthonormal bases in the three subspaces, each such tensor ${\cal Y}$ can be written uniquely as ${\cal Y}=({\bf S},{\bf T},{\bf U})\cdot{\cal G}$, with ${\cal G}$ an $R_1\times R_2\times R_3$ tensor. This provides the homeomorphism to Euclidian space. By a general result of \cite{Dra} over the complex field, for almost all ${\cal Z}$ the number of (complex) stationary points in the approximation from $M_{(R_1,R_2,R_3)}(I,J,K)$ with multilinear rank exactly $(R_1,R_2,R_3)$ is finite, and the real stationary points are a subset of the complex stationary points. It remains to consider stationary points with multilinear rank less than $(R_1,R_2,R_3)$. For each triplet $(r_1,r_2,r_3)<(R_1,R_2,R_3)$ and set $M_{(r_1,r_2,r_3)}(I,J,K)$ we again apply the result of \cite{Dra} to obtain a finite number of stationary points with multilinear rank exactly $(r_1,r_2,r_3)$. Note that the stationary points in the approximation from $M_{(r_1,r_2,r_3)}(I,J,K)$ may not be stationary points in the approximation from $M_{(R_1,R_2,R_3)}(I,J,K)$, whereas the converse is true for the stationary points that have multilinear rank at most $(r_1,r_2,r_3)$. Since the number of triplets $(r_1,r_2,r_3)$ is finite and we obtain finitely many stationary points for each triplet, it follows that there are finitely many stationary points in the approximation from $M_{(R_1,R_2,R_3)}(I,J,K)$. This completes the proof of Lemma~\ref{lem-M-AG}.

For Lemma~\ref{lem-S-AG} we need to consider the smooth points of $\ol{S}_2(I,J,K)$. The tensors in $\ol{S}_2(I,J,K)$ are either of the form ${\bf a}_1\circ{\bf b}_1\circ{\bf c}_1+{\bf a}_2\circ{\bf b}_2\circ{\bf c}_2$ or of the form ${\bf a}_1\circ{\bf b}_1\circ{\bf c}_2+{\bf a}_1\circ{\bf b}_2\circ{\bf c}_1+{\bf a}_2\circ{\bf b}_1\circ{\bf c}_1$ \cite{DSL}. The first type are smooth points only when each pair of ${\bf a}_i$, ${\bf b}_i$, and ${\bf c}_i$ has rank 2. That is, when the tensor has multilinear rank (2,2,2). The second type have rank 3 if and only if each pair of ${\bf a}_i$, ${\bf b}_i$, and ${\bf c}_i$ has rank 2 \cite{DSL}. However, the vectors ${\bf a}_i$, ${\bf b}_i$, and ${\bf c}_i$ are not unique \cite[lemma 4.1]{Ste-TV}. Hence, these are singular points in $\ol{S}_2(I,J,K)$. It is well known that the singular points Sing$(V)$ of an irreducible variety $V$ form a strict subvariety of $V$ that is the union of finite irreducible varieties $V_1,\ldots,V_m$. By the result of \cite{Dra} the number of stationary points in the approximation from $V\backslash$ Sing$(V)$ is finite for almost all ${\cal Z}$. Next, we consider the approximation from $V_j$ and obtain that the number of stationary points in $V_j\backslash$ Sing$(V_j)$ is finite by \cite{Dra}, for $j=1,\ldots,m$, and almost all ${\cal Z}$. We can continue this process and consider the approximation from Sing$(V_j)$, and so forth. Note that a stationary point in the approximation from $V_j$ may not be a stationary point in the approximation from $V$. Since $\ol{S}_2(I,J,K)$ is an irreducible variety, we obtain that for almost all ${\cal Z}$ the number of stationary points in the approximation from $\ol{S}_2(I,J,K)$ is finite.

\section*{Appendix B: proof of Lemma~\ref{lem-S-rank222}}
\refstepcounter{section}
\setcounter{equation}{0}
\renewcommand{\thesection}{B}
\renewcommand{\theequation}{B.\arabic{equation}}
For columnwise orthogonal ${\bf S}\in\R^{I\times 2}$, ${\bf T}\in\R^{J\times 2}$, and ${\bf U}\in\R^{K\times 2}$, and ${\cal G}\in\R^{2\times 2\times 2}$, let ${\cal X}=({\bf S},{\bf T},{\bf U})\cdot{\cal G}$ be a local minimum in the approximation of ${\cal Z}\in\R^{I\times J\times K}$ from $\ol{S}_2(I,J,K)$, with ${\cal Z}\notin\ol{S}_2(I,J,K)$. Then rank$({\cal X})\ge 2$ is easy to show using the proof of Lemma~\ref{lem-optrank}. From the orbits of $\R^{2\times 2\times 2}$ in Appendix D it follows that mrank$({\cal X})=$ mrank$({\cal G})$ equals either (1,2,2), (2,1,2), (2,2,1), or (2,2,2). Below, we prove that mrank$({\cal X})=(2,2,1)$ does not occur for almost all ${\cal Z}$. The proofs for (1,2,2) and (2,1,2) are analogous. We prove the following results.

\begin{lem}
\label{lem-S-rank221}
Let ${\cal Z}\in\R^{I\times J\times K}$, with ${\cal Z}\notin\ol{S}_2(I,J,K)$ and {\rm mrank}$({\cal Z})\ge (2,2,2)$. Let ${\cal X}$ be a local minimum in the approximation of ${\cal Z}$ from $\ol{S}_2(I,J,K)$ with {\rm mrank}$({\cal X})=(2,2,1)$. Then orthogonal $\wt{\bf S}\in\R^{I\times I}$, $\wt{\bf T}\in\R^{J\times J}$, and $\wt{\bf U}\in\R^{K\times K}$ exist such that $\wt{\cal Z}=(\wt{\bf S}^T,\wt{\bf T}^T,\wt{\bf U}^T)\cdot{\cal Z}$ has $I\times J$ slices $\wt{\bf Z}_k$, with $\wt{\bf Z}_k=\left[\begin{array}{cc} {\bf O} & {\bf O}\\ {\bf O} & {\bf M}_k\end{array}\right]$ for $k=2,\ldots,K$, where ${\bf M}_k\in\R^{(I-2)\times (J-2)}$. Furthermore, $(\wt{\bf Z}_1)_{ij}=0$ for $i\neq j$.
\ep
\end{lem}

\begin{lem}
\label{lem-S-formZ}
For almost all ${\cal Z}\in\R^{I\times J\times K}$ with ${\cal Z}\notin\ol{S}_2(I,J,K)$, no orthogonal $\wt{\bf S}\in\R^{I\times I}$, $\wt{\bf T}\in\R^{J\times J}$, and $\wt{\bf U}\in\R^{K\times K}$ exist such that $\wt{\cal Z}=(\wt{\bf S}^T,\wt{\bf T}^T,\wt{\bf U}^T)\cdot{\cal Z}$ has the form in the statement of {\rm Lemma~\ref{lem-S-rank221}}.
\ep
\end{lem}

\noindent We start by proving the following auxiliary result.
\begin{lem}
\label{lem-localrank}
For ${\bf Z}\in\R^{I\times J}$ with {\rm rank}$({\bf Z})\ge R\ge 1$, let ${\bf X}$ be a locally best {\rm rank}-$R$ approximation. Then ${\bf X}$ is a $($globally$)$ best {\rm rank}-$R$ approximation of ${\bf Z}$.
\end{lem}

\noindent {\bf Proof.} The problem of finding a best rank-$R$ approximation of ${\bf Z}$ is identical to finding a best rank-$(R,R,1)$ approximation of ${\bf Z}$ as described in section 2. That is, we are looking for orthogonal $\wt{\bf S}\in\R^{I\times I}$ and $\wt{\bf T}\in\R^{J\times J}$ such that $\wt{\bf Z}=\wt{\bf S}^T\,{\bf Z}\,\wt{\bf T}=\left[\begin{array}{cc} {\bf G} & {\bf L}\\ {\bf N} & {\bf M}\end{array}\right]$, with ${\bf G}\in\R^{R\times R}$, ${\bf L}\in\R^{R\times (J-R)}$, ${\bf N}\in\R^{(I-R)\times R}$, and ${\bf M}\in\R^{(I-R)\times (J-R)}$, has maximal $||{\bf G}||^2_F$. The corresponding best rank-$R$ approximation is then given by ${\bf X}={\bf S}\,{\bf G}\,{\bf T}^T$, where ${\bf S}$ and ${\bf T}$ consist of the first $R$ columns of $\wt{\bf S}$ and $\wt{\bf T}$, respectively. 

Let ${\bf X}={\bf S}\,{\bf G}\,{\bf T}^T$ be a locally best rank-$R$ approximation of ${\bf Z}$. As in the proof of Lemma~\ref{lem-optrank}, it can be seen that rank$({\bf X})=$ rank$({\bf G})=R$. First-order conditions (\ref{eq-statM-row}) and (\ref{eq-statM-col}) imply that ${\bf N}={\bf O}$ and ${\bf L}={\bf O}$, respectively. Using the SVDs of ${\bf G}$ and ${\bf M}$, we may assume without loss of generality that ${\bf G}={\rm diag}(\sigma_1({\bf G}),\ldots,\sigma_R({\bf G}))$, with $\sigma_1({\bf G})\ge\ldots\ge\sigma_R({\bf G})>0$, and ${\bf M}=\left[\begin{array}{cc} {\rm diag}(\sigma_1({\bf M}),\ldots,\sigma_{R_z-R}({\bf M})) & {\bf O} \\ {\bf O} & {\bf O}\end{array}\right]$, with $\sigma_1({\bf M})\ge\ldots\ge\sigma_{R_z-R}({\bf M})>0$, and $R_z={\rm rank}({\bf Z})$. Hence, the singular values of ${\bf Z}$ are given by $\sigma_1({\bf G}),\ldots,$\\ $\sigma_R({\bf G}),\sigma_1({\bf M}),\ldots,\sigma_{R_z-R}({\bf M})$. By \cite{EY} we need to show that $\sigma_R({\bf G})\ge\sigma_1({\bf M})$. For $R_z=R$ we have ${\bf X}={\bf Z}$ and we are done. In the following, let $R_z>R$. Let ${\bf e}_{i,n}\in\R^n$ be the $i$th unit vector. For small $t$, consider the rank-$R$ matrix 
\be
{\bf X}(t)=\wt{\bf S}\,\left[\left(\sum_{r=1}^{R-1}\sigma_r({\bf G})\,{\bf e}_{r,I}\,{\bf e}_{r,J}^T\right)+\sigma_R({\bf G})\,({\bf e}_{R,I}+t\,{\bf e}_{R+1,I})\,({\bf e}_{R,J}+t\,{\bf e}_{R+1,J})^T\right]\,\wt{\bf T}^T\,.
\ee

\noindent We have
\begin{eqnarray}
\label{eq-normdiff}
||{\bf Z}-{\bf X}(t)||^2_F &=& \sum_{r=2}^{R_z-R}\sigma_r^2({\bf M})+(\sigma_1({\bf M})-t^2\,\sigma_R({\bf G}))^2+2\,t^2\,\sigma_R^2({\bf G}) \nonumber \\
&=& ||{\bf Z}-{\bf X}||^2_F +t^2\,\sigma_R({\bf G})\,(t^2\,\sigma_R({\bf G})+2\,(\sigma_R({\bf G})-\sigma_1({\bf M})))\,.
\end{eqnarray}

\noindent By local optimality of ${\bf X}$ we must have $||{\bf Z}-{\bf X}(t)||^2_F\ge||{\bf Z}-{\bf X}||^2_F$ for small enough $t$. Hence, $t^2\,\sigma_R({\bf G})+2\,(\sigma_R({\bf G})-\sigma_1({\bf M}))\ge 0$ for small enough $t$, which implies $\sigma_R({\bf G})\ge\sigma_1({\bf M})$. This completes the proof.
\ep

\noindent {\bf Proof of Lemma~\ref{lem-S-rank221}}\\
\noindent From the orbits of $\R^{2\times 2\times 2}$ in Appendix D it follows that mrank$({\cal X})=(2,2,1)$ implies rank$({\cal X})=2$. An orthogonal $\wt{\bf U}\in\R^{K\times K}$ exists such that $\wh{\cal X}=({\bf I}_I,{\bf I}_J,\wt{\bf U}^T)\cdot{\cal X}$ has $I\times J$ slices $\wh{\bf X}_k$ with rank$(\wh{\bf X}_1)=2$ and $\wh{\bf X}_k={\bf O}$ for $k=2,\ldots,K$. Let $\wh{\cal Z}=({\bf I}_I,{\bf I}_J,\wt{\bf U}^T)\cdot{\cal Z}$ have $I\times J$ slices $\wh{\bf Z}_k$, for $k=1,\ldots,K$. Then $\wh{\cal X}$ is a local minimum in the approximation of $\wh{\cal Z}$ from $\ol{S}_2(I,J,K)$ with {\rm mrank}$(\wh{\cal X})=(2,2,1)$. We have $||\wh{\cal Z}-\wh{\cal X}||^2_F=||\wh{\bf Z}_1-\wh{\bf X}_1||^2_F+\sum_{k=2}^K ||\wh{\bf Z}_k||^2_F$. Local optimality of $\wh{\cal X}$ implies that $\wh{\bf X}_1$ is a locally best rank-2 approximation of $\wh{\bf Z}_1$. Indeed, otherwise a locally better approximation from $\ol{S}_2(I,J,K)$ can be found by varying $\wh{\bf X}_1$ in its neighborhood of rank-2 matrices. Let the SVD of $\wh{\bf Z}_1$ be given by $\wh{\bf Z}_1=\wt{\bf S}\,{\rm diag}(\sigma_1(\wh{\bf Z}_1),\sigma_2(\wh{\bf Z}_1),\ldots,\sigma_{\min(I,J)}(\wh{\bf Z}_1))\,\wt{\bf T}^T$, where the singular values $\sigma_i(\wh{\bf Z}_1)$ are in descending order, and $\wt{\bf S}\in\R^{I\times I}$ and $\wt{\bf T}\in\R^{J\times J}$ are orthogonal. It follows from Lemma~\ref{lem-localrank} and \cite{EY} that $\wh{\bf X}_1$ is of the form $\wh{\bf X}_1=\wt{\bf S}\,{\rm diag}(\sigma_{1}(\wh{\bf Z}_1),\sigma_{2}(\wh{\bf Z}_1),0,\ldots,0)\,\wt{\bf T}^T$. Note that equal nonzero singular values imply nonuniqueness of the corresponding left- and right singular vectors in $\wt{\bf S}$ and $\wt{\bf T}$. However, since this nonuniqueness is present in both the SVD of $\wh{\bf Z}_1$ and the SVD of a best rank-2 approximation $\wh{\bf X}_1$ when one or both of $\sigma_{1}(\wh{\bf Z}_1)$ and $\sigma_{2}(\wh{\bf Z}_1)$ have multiplicity greater than one, it remains true that $\wh{\bf X}_1$ can be written in the form above. Now define $\wt{\cal Z}=(\wt{\bf S}^T,\wt{\bf T}^T,{\bf I}_K)\cdot\wh{\cal Z}$ and $\wt{\cal X}=(\wt{\bf S}^T,\wt{\bf T}^T,{\bf I}_K)\cdot\wh{\cal X}$. This proves the form of $\wt{\bf Z}_1=\wt{\bf S}^T\,\wh{\bf Z}_1\,\wt{\bf T}$, and $\wt{\cal X}$ is a local minimum in the approximation of $\wt{\cal Z}$ from $\ol{S}_2(I,J,K)$.

Let ${\bf e}_{i,n}\in\R^n$ be the $i$th unit vector, and set $g_{111}=\sigma_{1}(\wh{\bf Z}_1)$ and $g_{221}=\sigma_{2}(\wh{\bf Z}_1)$ for ease of presentation. 
Then $\wt{\bf X}_1=g_{111}\,{\bf e}_{1,I}\,{\bf e}_{1,J}^T+g_{221}\,{\bf e}_{2,I}\,{\bf e}_{2,J}^T$. Also, we have $\wt{\bf X}_k={\bf O}$ for $k=2,\ldots,K$. To prove the form of $\wt{\bf Z}_k$, $k=2,\ldots,K$, we consider small perturbations of $\wt{\cal X}$ in $\ol{S}_2(I,J,K)$ and use the local optimality of $\wt{\cal X}$. Let $k\in\{2,\ldots,K\}$ be fixed. For parameters $s,t\in\R$, the perturbation of $\wt{\bf X}_1$ is given by
\be
\label{eq-X1per}
\wt{\bf X}_1(s,t)=g_{111}\,({\bf e}_{1,I}+s\,{\bf h}_{1,1})\,({\bf e}_{1,J}+t\,{\bf h}_{1,2})^T+g_{221}\,({\bf e}_{2,I}+s\,{\bf h}_{2,1})\,({\bf e}_{2,J}+t\,{\bf h}_{2,2})^T\,,
\ee

\noindent with ${\bf h}_{i,1}\in$ span$({\bf e}_{3,I},\ldots,{\bf e}_{I,I})$ and ${\bf h}_{i,2}\in$ span$({\bf e}_{3,J},\ldots,{\bf e}_{J,J})$, for $i=1,2$. For parameter $u\in\R$, the perturbation of $\wt{\bf X}_k$ is given by $\wt{\bf X}_k(s,t,u)=u\,{\bf a}(s)\,{\bf b}(t)^T$, with ${\bf a}(s)\in$ span$({\bf e}_{1,I}+s\,{\bf h}_{1,1},\,{\bf e}_{2,I}+s\,{\bf h}_{2,1})$ and ${\bf b}(t)\in$ span$({\bf e}_{1,J}+t\,{\bf h}_{1,2},\,{\bf e}_{2,J}+t\,{\bf h}_{2,2})$. We do not perturb $\wt{\bf X}_m={\bf O}$ for $m\neq k$ and $m\ge 2$. We claim that the perturbation $\wt{\cal X}(s,t,u)$ lies in $\ol{S}_2(I,J,K)$ for small $s,t,u$. It suffices to show that the $I\times J\times 2$ tensor with slices $\wt{\bf X}_1(s,t)$ and $\wt{\bf X}_k(s,t,u)$ lies in $\ol{S}_2(I,J,2)$. In $\R^I$, $\R^J$, and $\R^2$ we choose the following bases: ${\bf e}_{1,I}+s\,{\bf h}_{1,1},\,{\bf e}_{2,I}+s\,{\bf h}_{2,1},\,{\bf e}_{3,I},\ldots,{\bf e}_{I,I}$ in $\R^I$, ${\bf e}_{1,J}+t\,{\bf h}_{1,2},\,{\bf e}_{2,J}+t\,{\bf h}_{2,2},\,{\bf e}_{3,J},\ldots,{\bf e}_{J,J}$ in $\R^J$, and ${\bf e}_{1,2},\,{\bf e}_{2,2}$ in $\R^2$. With respect to these bases, the $I\times J\times 2$ tensor has representation with $I\times J$ slices 
\be
\label{eq-reprIJ2}
\left[\begin{array}{cc} {\rm diag}(g_{111},\,g_{221}) & {\bf O}\\ {\bf O} & {\bf O}\end{array}\right]\,,\quad\quad
\left[\begin{array}{cc} u\,{\bf a}(s)\,{\bf b}(t)^T & {\bf O}\\ {\bf O} & {\bf O}\end{array}\right]\,.
\ee

\noindent Hence, we need to show that the $2\times 2\times 2$ tensor with slices equal to the nonzero parts in (\ref{eq-reprIJ2}) lies in $\ol{S}_2(2,2,2)$. But this follows from Proposition~\ref{p-1} since slice 2 multiplied by the inverse of slice 1 cannot have complex eigenvalues. 

Next, we use the perturbation $\wt{\cal X}(s,t,u)$ and local optimality of $\wt{\cal X}$ to prove the form of $\wt{\bf Z}_k$. We have by the triangle inequality that
\begin{eqnarray}
\label{eq-normdecomp}
||\wt{\cal Z}-\wt{\cal X}(s,t,u)||^2_F &=& ||\wt{\bf Z}_1-\wt{\bf X}_1(s,t)||^2_F+
\sum_{m\neq k,\,m\ge 2}||\wt{\bf Z}_m-\wt{\bf X}_m||^2_F +
||\wt{\bf Z}_k-\wt{\bf X}_k(s,t,u)||^2_F \nonumber \\
 &\le& \sum_{m\neq k}||\wt{\bf Z}_m-\wt{\bf X}_m||^2_F+||\wt{\bf X}_1-\wt{\bf X}_1(s,t)||^2_F+||\wt{\bf Z}_k-\wt{\bf X}_k(s,t,u)||^2_F\,.
\end{eqnarray}

\noindent For small enough $s,t,u$, local optimality of $\wt{\cal X}$ thus implies $||\wt{\bf X}_1-\wt{\bf X}_1(s,t)||^2_F+||\wt{\bf Z}_k-\wt{\bf X}_k(s,t,u)||^2_F\ge ||\wt{\bf Z}_k-\wt{\bf X}_k||^2_F=||\wt{\bf Z}_k||^2_F$. 

Let $\wt{\bf Z}_k=\left[\begin{array}{cc} {\bf G}_k & {\bf L}_k\\ {\bf N}_k & {\bf M}_k\end{array}\right]$, with ${\bf G}_k\in\R^{2\times 2}$, ${\bf L}_k\in\R^{2\times (J-2)}$, ${\bf N}_k\in\R^{(I-2)\times 2}$, and ${\bf M}_k\in\R^{(I-2)\times (J-2)}$. To prove Lemma~\ref{lem-S-rank221}, we need to show that ${\bf G}_k={\bf O}$, ${\bf L}_k={\bf O}$, and ${\bf N}_k={\bf O}$. Set $s=t=0$. Then choosing ${\bf a}(0){\bf b}(0)^T={\bf e}_{i,I},{\bf e}_{j,J}^T$ for $i,j\in\{1,2\}$ implies $||\wt{\bf Z}_k-\wt{\bf X}_k(0,0,u)||^2_F<||\wt{\bf Z}_k||^2_F$ for small enough $u$ when $({\bf G}_k)_{ij}\neq 0$. Hence, we obtain ${\bf G}_k={\bf O}$. We write ${\bf a}(s)=\left(\begin{array}{c} {\bf a}_1\\ s\,{\bf a}_2\end{array}\right)$ and ${\bf b}(t)=\left(\begin{array}{c} {\bf b}_1\\ t\,{\bf b}_2\end{array}\right)$, with ${\bf a}_1,{\bf b}_1\in$ span$({\bf e}_{1,2},{\bf e}_{2,2})$, ${\bf a}_2\in$ span$({\bf e}_{1,I-2},\ldots,{\bf e}_{I-2,I-2})$, and ${\bf b}_2\in$ span$({\bf e}_{1,J-2},\ldots,{\bf e}_{J-2,J-2})$. It follows that 
\be
\label{eq-Xkper}
\wt{\bf X}_k(s,t,u)=u\,{\bf a}(s)\,{\bf b}(t)^T=\left[\begin{array}{cc}
u\,{\bf a}_1\,{\bf b}_1^T & t\,u\,{\bf a}_1\,{\bf b}_2^T\\
s\,u\,{\bf a}_2\,{\bf b}_1^T & s\,t\,u\,{\bf a}_2\,{\bf b}_2^T\end{array}\right]\,,
\ee

\noindent and $||\wt{\bf Z}_k-\wt{\bf X}_k(s,t,u)||^2_F$ can be written as
\be
\label{eq-normk}
||u\,{\bf a}_1\,{\bf b}_1^T||^2_F+||{\bf L}_k-t\,u\,{\bf a}_1\,{\bf b}_2^T||^2_F+
||{\bf N}_k-s\,u\,{\bf a}_2\,{\bf b}_1^T||^2_F+
||{\bf M}_k-s\,t\,u\,{\bf a}_2\,{\bf b}_2^T||^2_F\,.
\ee

\noindent Let $s=0$ and fix $t,u\neq 0$. Then (\ref{eq-normk}) equals $||{\bf L}_k||^2_F+||{\bf N}_k||^2_F+||{\bf M}_k||^2_F-2\,t\,u\,{\bf a}_1^T{\bf L}_k\,{\bf b}_2+t^2\,u^2\,({\bf a}_1^T{\bf a}_1)\,({\bf b}_2^T{\bf b}_2)+u^2\,({\bf a}_1^T{\bf a}_1)\,({\bf b}_1^T{\bf b}_1)$. We have $||\wt{\bf X}_1-\wt{\bf X}_1(0,t)||^2_F=t^2\,(g_{111}^2\,{\bf h}_{1,2}^T{\bf h}_{1,2}+g_{221}^2\,{\bf h}_{2,2}^T{\bf h}_{2,2})$. Local optimality of $\wt{\cal X}$ implies that 
\bdm
\left(\frac{t}{u}\right)^2\,(g_{111}^2\,{\bf h}_{1,2}^T{\bf h}_{1,2}+g_{221}^2\,{\bf h}_{2,2}^T{\bf h}_{2,2})-2\,\left(\frac{t}{u}\right)\,{\bf a}_1^T{\bf L}_k\,{\bf b}_2+({\bf a}_1^T{\bf a}_1)\,({\bf b}_1^T{\bf b}_1)
\edm
\be
\label{eq-optcond}
+\;t^2\,({\bf a}_1^T{\bf a}_1)\,({\bf b}_2^T{\bf b}_2)\ge 0\,,
\ee

\noindent for $t,u$ small enough. Hence, the first three terms of (\ref{eq-optcond}) together should be nonnegative for $t,u$ small enough. Denote this second degree polynomial by $c_2\,x^2-2\,c_1\,x+c_0$, where $x=t/u$ can have any value. We have $c_2>0$ and $c_0>0$ and the minimal value of the polynomial is $c_0-c_1^2/c_2$ at $x=c_1/c_2$. Hence, the polynomial is nonnegative when $c_0\ge c_1^2/c_2$. For any ${\bf a}_1^T{\bf L}_k\neq {\bf 0}^T$ we can choose ${\bf b}_2$ such that $c_1^2>c_0\,c_2$. Thus ${\bf a}_1^T\,{\bf L}_k={\bf 0}^T$ for any ${\bf a}_1\in$ span$({\bf e}_{1,2},{\bf e}_{2,2})$. This is only possible when ${\bf L}_k={\bf O}$. Analogously, it can be proven that ${\bf N}_k={\bf O}$. This shows that $\wt{\bf Z}_k$ is of the form $\left[\begin{array}{cc} {\bf O} & {\bf O}\\ {\bf O} & {\bf M}_k\end{array}\right]$. Since $k\in\{2,\ldots,K\}$ is arbitrary, this completes the proof.

Note that the construction of the perturbation $\wt{\cal X}(s,t,u)$ requires $\min(I,J)\ge 3$. When $\min(I,J)=2$ Lemma~\ref{lem-S-rank221} states that $\wt{\bf Z}_k={\bf O}$ for $k=2,\ldots,K$. This is proven as follows. When $I=2$ and $J\ge 3$ we set $s=0$ and ${\bf N}_k$ and ${\bf M}_k$ do not exist, and ${\bf G}_k={\bf O}$ and ${\bf L}_k={\bf O}$ follow as above. When $J=2$ and $I\ge 3$ the result is obtained analogously by setting $t=0$. When $I=J=2$ we set $s=t=0$ and $\wt{\bf Z}_k={\bf G}_k={\bf O}$ follows as above. Note that the proof of the diagonal form of $\wt{\bf Z}_1$ is still valid when $\min(I,J)=2$.
\ep

\noindent {\bf Proof of Lemma~\ref{lem-S-formZ}}\\
\noindent Let $\wt{\bf S}\in\R^{I\times I}$, $\wt{\bf T}\in\R^{J\times J}$, and $\wt{\bf U}\in\R^{K\times K}$ be orthogonal such that $\wt{\cal Z}=(\wt{\bf S}^T,\wt{\bf T}^T,\wt{\bf U}^T)\cdot{\cal Z}$ has the form in the statement of Lemma~\ref{lem-S-rank221}. When $\min(I,J)=2$ we have $\wt{\bf Z}_k={\bf O}$ for $k=2,\ldots,K$ and, hence, the mode-3 rank of $\wt{\cal Z}$ equals 1. Since this is equal to the mode-3 rank of ${\cal Z}$, it does not occur for almost all ${\cal Z}$. In the following, let $\min(I,J)\ge 3$. 

Let $\wt{\bf U}^T=\left[\begin{array}{cc} \tilde{u}_{11} & \tilde{\bf u}_{21}^T\\ \tilde{\bf u}_{12} & \wt{\bf U}_{22}\end{array}\right]$, with $\tilde{\bf u}_{12},\tilde{\bf u}_{21}\in\R^{K-1}$, and $\wt{\bf U}_{22}\in\R^{(K-1)\times (K-1)}$. We have
\be
\label{eq-unfoldZ}
[\wt{\bf Z}_2\;\ldots\;\wt{\bf Z}_K]=\wt{\bf S}^T\,[{\bf Z}_1\;\ldots\;{\bf Z}_K]\,([\tilde{\bf u}_{12}\;\wt{\bf U}_{22}]^T\otimes \wt{\bf T})\,,
\ee

\noindent where $\otimes$ denotes the Kronecker product. By assumption, the matrix on the left-hand side of (\ref{eq-unfoldZ}) has allzero rows $1$ and $2$ and, hence, rank at most $I-2$. Let ${\bf L}$ be the inverse of a nonsingular $(K-1)\times (K-1)$ submatrix of $[\tilde{\bf u}_{12}\;\wt{\bf U}_{22}]^T$, which has rank $K-1$. Then 
\be
\label{eq-unfoldZ2}
\wt{\bf S}\,[\wt{\bf Z}_2\;\ldots\;\wt{\bf Z}_K]\,({\bf L}\otimes\wt{\bf T}^T)=
[{\bf Z}_1\;\ldots\;{\bf Z}_K]\,(([\boldsymbol{\beta}\;{\bf I}_{K-1}]\,{\bf\Pi})^T\otimes {\bf I}_J)\,,
\ee

\noindent where ${\bf\Pi}\in\R^{K\times K}$ is a permutation matrix and $\boldsymbol{\beta}\in\R^{K-1}$. Since the transformations are nonsingular, the matrix in (\ref{eq-unfoldZ2}) also has rank at most $I-2$. Without loss of generality we set ${\bf\Pi}={\bf I}_K$. Then (\ref{eq-unfoldZ2}) equals
\be
\label{eq-unfoldZ3}
[{\bf Z}_1\;\ldots\;{\bf Z}_K]\,([\boldsymbol{\beta}\;{\bf I}_{K-1}]^T\otimes {\bf I}_J)=
[{\bf Z}_2+\beta_1\,{\bf Z}_1\;\ldots\;{\bf Z}_K+\beta_{K-1}\,{\bf Z}_1]\,.
\ee

\noindent Hence, for some $[{\bf v}\;{\bf w}]$ with rank$([{\bf v}\;{\bf w}])=2$ we have $[{\bf v}\;{\bf w}]^T\,({\bf Z}_k+\beta_{k-1}\,{\bf Z}_1)={\bf O}$ for $k=2,\ldots,K$. When $I=J$ this implies that $-\beta_{k-1}$ is a real eigenvalue of $({\bf Z}_1^T)^{-1}{\bf Z}_k^T$ with associated eigenvectors ${\bf v}$ and ${\bf w}$. However, for almost all ${\cal Z}$, the matrix $({\bf Z}_1^T)^{-1}{\bf Z}_k^T$ has $I=J$ distinct eigenvalues with one associated eigenvector each. This completes the proof for $I=J$. For $I<J$ the same arguments can be used for the submatrices of ${\bf Z}_k$ and ${\bf Z}_1$ consisting of the first $I$ columns.

Since $\wt{\cal Z}=(\wt{\bf S}^T,\wt{\bf T}^T,\wt{\bf U}^T)\cdot{\cal Z}$ has the form in the statement of Lemma~\ref{lem-S-rank221}, the matrix $[\wt{\bf Z}_2^T\;\ldots\;\wt{\bf Z}_K^T]$ also has allzero rows $1$ and $2$. Analogous to (\ref{eq-unfoldZ})--(\ref{eq-unfoldZ2}) we obtain that $[{\bf Z}_1^T\;\ldots\;{\bf Z}_K^T]\,(([\boldsymbol{\alpha}\;{\bf I}_{K-1}]\,\wh{\bf\Pi})^T\otimes {\bf I}_I)$ has rank at most $J-2$ for some $\boldsymbol{\alpha}\in\R^{K-1}$ and permutation matrix $\wh{\bf\Pi}\in\R^{K\times K}$. Analogous to the above, this is not possible for almost all ${\cal Z}$ when $J\le I$. This completes the proof.
\ep

\section*{Appendix C: proof of Lemma~\ref{lem-M-interior}}
\refstepcounter{section}
\setcounter{equation}{0}
\renewcommand{\thesection}{C}
\renewcommand{\theequation}{C.\arabic{equation}}
For columnwise orthogonal ${\bf S}\in\R^{I\times 2}$, ${\bf T}\in\R^{J\times 2}$, and ${\bf U}\in\R^{K\times 2}$, and ${\cal G}\in\R^{2\times 2\times 2}$, let ${\cal X}=({\bf S},{\bf T},{\bf U})\cdot{\cal G}$ be a locally best rank-(2,2,2) approximation of ${\cal Z}\in\R^{I\times J\times K}$, with mrank$({\cal Z})>(2,2,2)$. Then rank$({\cal X})\ge 2$ is easy to show using the proof of Lemma~\ref{lem-optrank}. From the orbits of $\R^{2\times 2\times 2}$ in Appendix D it follows that mrank$({\cal X})=$ mrank$({\cal G})$ equals either (1,2,2), (2,1,2), (2,2,1), or (2,2,2). Lemma~\ref{lem-S-rank222} implies that mrank$({\cal X})<(2,2,2)$ does not occur for almost all ${\cal Z}$. Indeed, mrank$({\cal X})<(2,2,2)$ and rank$({\cal X})\ge 2$ imply rank$({\cal X})=2$. In that case ${\cal X}$ is also a local minimum in the best approximation from $S_2(I,J,K)\subset M_{(2,2,2)}(I,J,K)$ and Lemma~\ref{lem-S-rank222} applies. Hence, mrank$({\cal X})=(2,2,2)$ and rank$({\cal X})$ equals 2 or 3. Recall from Appendix A that this implies that ${\cal X}$ is a smooth point of the variety $M_{(2,2,2)}(I,J,K)$. We need to show that for almost all ${\cal Z}$ the local minimizer ${\cal X}$ does not lie on the boundary of $\ol{S}_2(I,J,K)$. The latter is a subvariety of $M_{(2,2,2)}(I,J,K)$ of codimension 1 (see Appendix D). Recall from Appendix A that the number of smooth stationary points in $M_{(2,2,2)}(I,J,K))$ is finite.

Denote the tangent space of $M_{(2,2,2)}(I,J,K)$ at ${\cal X}$ by $T_{\cal X}\subset\R^{I\times J\times K}$. We have $T_{\cal X}=\{{\cal Y}\in\R^{I\times J\times K}:\;{\cal Y}={\cal X}+{\cal X}_2,\;{\cal X}_2\in W\}$ for some subspace $W\subset M_{(2,2,2)}(I,J,K))$. Since ${\cal X}$ is a smooth point, it follows that dim$(T_{\cal X})=$ dim$(W)=$ dim$(M_{(2,2,2)}(I,J,K))$. Since ${\cal X}$ is a smooth stationary point, the gradient of $||{\cal Z}-{\cal Y}||^2_F$ restricted to $T_{\cal X}$ is zero at ${\cal X}$. 
Generic ${\cal Z}\in\R^{I\times J\times K}$ with $(I,J,K)>(2,2,2)$ for which ${\cal X}$ is a stationary point are of the form ${\cal Z}={\cal X}+{\cal X}_3$ for some ${\cal X}_3\in W^{\perp}$, with $W^{\perp}$ denoting the orthogonal complement of $W$ in $\R^{I\times J\times K}$. Hence, the variety of such ${\cal Z}$ has dimension $IJK-{\rm dim}(W)=IJK-{\rm dim}(M_{(2,2,2)}(I,J,K))$. Since ${\cal X}$ is an isolated smooth stationary point, it follows from the above that by varying ${\cal Z}$ in a neighborhood in $\R^{I\times J\times K}$ we obtain a neighborhood of ${\cal X}$ in $M_{(2,2,2)}(I,J,K)$ of dimension dim$(M_{(2,2,2)}(I,J,K))$. This contradicts that ${\cal X}$ is a boundary point of $\ol{S}_2(I,J,K)$ for generic ${\cal Z}$ since then the neighborhood of ${\cal X}$ in $M_{(2,2,2)}(I,J,K)$ should have  dimension (at most) dim$(M_{(2,2,2)}(I,J,K))-1$. This completes the proof.
\ep

\section*{Appendix D: classification of $2\times 2\times 2$ tensors}
\refstepcounter{section}
\setcounter{equation}{0}
\renewcommand{\thesection}{D}
\renewcommand{\theequation}{D.\arabic{equation}}
In \cite{DSL} it is shown that all tensors in $\R^{2\times 2\times 2}$ can be transformed to eight canonical forms, i.e., for each ${\cal Y}\in\R^{2\times 2\times 2}$ nonsingular ${\bf S}\in\R^{2\times 2}$, ${\bf T}\in\R^{2\times 2}$, and ${\bf U}\in\R^{2\times 2}$ exist such that $({\bf S},{\bf T},{\bf U})\cdot{\cal Y}$ equals a canonical form. Moreover, each ${\cal Y}\in\R^{2\times 2\times 2}$ can be transformed to one canonical form only. This implies a classification of $\R^{2\times 2\times 2}$ into eight orbits, with rank and mrank fixed on each orbit. In Table~\ref{tab-orbits} the canonical forms are listed, together with values of rank, mrank, and the Lebesgue measure for each orbit.

As can be seen, a generic tensor in $\R^{2\times 2\times 2}$ is either in orbit $G_2$ or in orbit $G_3$. Hence, rank 2 and rank 3 both occur on sets of positive measure. It is shown in \cite{DSL} that ${\cal Z}\in\R^{2\times 2\times 2}$ with rank$({\cal Z})=3$ does not have  a best rank-2 approximation. It is stated in \cite{DSL} that there are seven orbits in $\C^{2\times 2\times 2}$, where $D_0$, $D_1$, $D_2$, $D_2'$, $D_2''$, and $D_3$ have the same canonical form, rank, mrank, and measure as for $\R^{2\times 2\times 2}$. The seventh orbit consists of rank-2, mrank-(2,2,2) tensors and is the only orbit with positive measure in $\C^{2\times 2\times 2}$. We state the following result of \cite{Ste-arxiv-GSD} for real tensors.

\newpage
\begin{table}[t]
\begin{center}
\begin{tabular}{ccccc}
\hline 
orbit & canonical form & rank & mrank & measure \\
\hline \noalign{\vskip 2mm} 
$D_0$ & $\left[\begin{array}{cc|cc} 0 & 0 & 0 & 0\\ 0 & 0 & 0 & 0\end{array}\right]$ & 0 & (0,0,0) & 0 \\ \noalign{\vskip 2mm} 
$D_1$ & $\left[\begin{array}{cc|cc} 1 & 0 & 0 & 0\\ 0 & 0 & 0 & 0\end{array}\right]$ & 1 & (1,1,1) & 0 \\   \noalign{\vskip 2mm} 
$D_2$ & $\left[\begin{array}{cc|cc} 1 & 0 & 0 & 0\\ 0 & 1 & 0 & 0\end{array}\right]$ & 2 & (2,2,1) & 0 \\  \noalign{\vskip 2mm} 
$D_2'$ & $\left[\begin{array}{cc|cc} 1 & 0 & 0 & 1\\ 0 & 0 & 0 & 0\end{array}\right]$ & 2 & (2,1,2) & 0 \\  \noalign{\vskip 2mm} 
$D_2''$ & $\left[\begin{array}{cc|cc} 1 & 0 & 0 & 0\\ 0 & 0 & 1 & 0\end{array}\right]$ & 2 & (1,2,2) & 0 \\  \noalign{\vskip 2mm} 
$G_2$ & $\left[\begin{array}{cc|cc} 1 & 0 & 0 & 0\\ 0 & 0 & 0 & 1\end{array}\right]$ & 2 & (2,2,2) & $>0$ \\ \noalign{\vskip 2mm} 
$D_3$ & $\left[\begin{array}{cc|cc} 1 & 0 & 0 & 1\\ 0 & 1 & 0 & 0\end{array}\right]$ & 3 & (2,2,2) & 0 \\ \noalign{\vskip 2mm} 
$G_3$ & $\left[\begin{array}{cc|cc} 1 & 0 & 0 & -1\\ 0 & 1 & 1 & 0\end{array}\right]$ & 3 & (2,2,2) & $>0$ \\  \noalign{\vskip 2mm} 
\hline
\end{tabular}
\end{center}
\caption{Orbits of $\R^{2\times 2\times 2}$ proven by \cite{DSL}.} 
\label{tab-orbits}
\end{table}

\clearpage

\begin{propos}
\label{p-1}
Let ${\cal Y}\in\R^{R\times R\times 2}$.
\begin{itemize}
\item[{\rm (a)}] If there exists a ${\bf U}\in\R^{2\times 2}$ nonsingular such that
${\cal X}=({\bf I}_R,{\bf I}_R,{\bf U})\cdot{\cal Y}$ has nonsingular
slice ${\bf X}_1$, then
\begin{itemize}
\item[{\rm (a1)}] ${\cal Y}$ is an interior point of $S_R(R,R,2)$
if ${\bf X}_2{\bf X}_1^{-1}$ has $R$ distinct real eigenvalues, and {\rm rank}$({\cal Y})=R$.
\item[{\rm (a2)}] ${\cal Y}$ is a boundary point of $S_R(R,R,2)$
if ${\bf X}_2{\bf X}_1^{-1}$ has $R$ real eigenvalues but not all distinct, with {\rm rank}$({\cal Y})=R$ if and only if ${\bf X}_2{\bf X}_1^{-1}$ has $R$ linearly independent eigenvectors.
\item[{\rm (a3)}] ${\cal Y}$ is an exterior point of $S_R(R,R,2)$
if ${\bf X}_2{\bf X}_1^{-1}$ has at least one pair of complex eigenvalues.
\end{itemize}
\item[{\rm (b)}] If there does not exist a ${\bf U}\in\R^{2\times 2}$ nonsingular
such that ${\cal X}=({\bf I}_R,{\bf I}_R,{\bf U})\cdot{\cal Y}$ has
nonsingular slice ${\bf X}_1$, then ${\cal Y}$ is a boundary
point of $S_R(R,R,2)$. \ep
\end{itemize}
\end{propos}

\noindent In parallel with Proposition~\ref{p-1} for $R=2$ it can be shown that the boundary of $\ol{S}_2(2,2,2)$ consists of orbits $D_0$, $D_1$, $D_2$, $D_2'$, $D_2''$, and $D_3$. Orbit $G_2$ forms the interior of $\ol{S}_2(2,2,2)$, while orbit $G_3$ is the exterior. It follows that the boundary of $\ol{S}_2(I,J,K)$ is formed by ${\cal X}=({\bf S},{\bf T},{\bf U})\cdot{\cal G}$ with columnwise orthogonal ${\bf S}\in\R^{I\times 2}$, ${\bf T}\in\R^{J\times 2}$, and ${\bf U}\in\R^{K\times 2}$, and ${\cal G}\in\R^{2\times 2\times 2}$ in orbit $D_0$, $D_1$, $D_2$, $D_2'$, $D_2''$, or $D_3$. Interior points of $\ol{S}_2(I,J,K)$ have corresponding ${\cal G}$ in orbit $G_2$. Exterior points have corresponding ${\cal G}$ in orbit $G_3$.

The boundary of $\ol{S}_2(2,2,2)$ is characterized by the hyperdeterminant being zero, where the latter is a polynomial in the entries of the $2\times 2\times 2$ tensor \cite{DSL}. Orbits $G_2$ and $G_3$ have positive and negative hyperdeterminant, respectively. It follows that the boundary of $\ol{S}_2(2,2,2)$ has codimension 1 in $\R^{2\times 2\times 2}$. This implies that the boundary of $\ol{S}_2(I,J,K)$ has codimension 1 in $M_{(2,2,2)}(I,J,K)$ and is a subvariety of $M_{(2,2,2)}(I,J,K)$.

\end{document}